\documentclass[12pt,a4paper]{article}
\usepackage{epsfig,latexsym,amsfonts,amssymb,amsmath,amscd,
graphics,theorem,epic}
\oddsidemargin36pt \evensidemargin74pt  \theoremstyle{plain}
\newtheorem{lemma}{Lemma}
\newtheorem{proposition}[lemma]{Proposition}
\newtheorem{remark}[lemma]{Remark}
\newtheorem{example}[lemma]{Example}
\newtheorem{theorem}{Theorem}

\newtheorem{corollary}[lemma]{Corollary}
{\theorembodyfont{\rmfamily} 
\font\ncsc=cmcsc10  \font\ntt=cmtt12
\begin{document}
\baselineskip=15pt
\newcommand{\pperp}{\hbox{$\perp\hskip-6pt\perp$}}
\newcommand{\ssim}{\hbox{$\hskip-2pt\sim$}}
\newcommand{\N}{{\mathbb N}}\newcommand{\Delp}{{\Pi}}
\newcommand{\Z}{{\mathbb Z}}
\newcommand{\R}{{\mathbb R}}
\newcommand{\SSS}{{\mathbb S}}
\newcommand{\CB}{{\mathbb B}}
\newcommand{\CS}{{\mathbb B_1}}
\newcommand{\mm}{{\mathfrak m}}
\newcommand{\C}{{\mathbb C}}
\newcommand{\Q}{{\mathbb Q}}\newcommand{\T}{{\mathbb T}}\newcommand{\K}{{\mathbb
K}}\newcommand{\D}{{\mathcal D}}
\newcommand{\PP}{{\mathbb P}}
\newcommand{\st}{{*}}\newcommand{\coker}{{\operatorname{coker}}}
\newcommand{\Prec}{{\operatorname{Prec}}}
\newcommand{\ev}{{\operatorname{ev}}}
\newcommand{\Id}{{\operatorname{Id}}}\newcommand{\irr}{{\operatorname{irr}}}
\newcommand{\Sym}{{\operatorname{Sym}}}
\newcommand{\fix}{{\operatorname{fix}}}\newcommand{\re}{{\operatorname{re}}}
\newcommand{\ima}{{\operatorname{im}}}
\newcommand{\nod}{{\operatorname{nod}}}\newcommand{\lcm}{{\operatorname{lcm}}}
\newcommand{\oeps}{{\overline\eps}}\newcommand{\Area}{{\operatorname{Area}}}
\newcommand{\codim}{{\operatorname{codim}}}\newcommand{\sol}{{\operatorname{sol}}}
\newcommand{\oDel}{{\widetilde\Del}}
\newcommand{\odd}{{\operatorname{odd}}}
\newcommand{\real}{{\operatorname{Re}}}\newcommand{\ind}{{\operatorname{ind}}}
\newcommand{\conv}{{\operatorname{conv}}}
\newcommand{\Span}{{\operatorname{Span}}}
\newcommand{\Ker}{{\operatorname{Ker}}}
\newcommand{\Fix}{{\operatorname{Fix}}}\newcommand{\ord}{{\operatorname{ord}}}
\newcommand{\sign}{{\operatorname{sign}}}
\newcommand{\Log}{{\operatorname{Log}}}
\newcommand{\jet}{{\operatorname{jet}}}
\newcommand{\oi}{{\overline i}}
\newcommand{\oj}{{\overline j}}
\newcommand{\ob}{{\overline b}}
\newcommand{\os}{{\overline s}}
\newcommand{\oa}{{\overline a}}
\newcommand{\oy}{{\overline y}}
\newcommand{\ow}{{\overline w}}
\newcommand{\ou}{{\overline u}}
\newcommand{\ot}{{\overline t}}
\newcommand{\oz}{{\overline z}}
\newcommand{\cheb}{\ensuremath{\mathrm{cheb}}}
\newcommand{\newi}{i}\newcommand{\bd}{{\boldsymbol d}}
\newcommand{\newj}{j}
\newcommand{\newm}{m}
\newcommand{\newl}{{\ell}}
\newcommand{\bw}{{\boldsymbol w}}\newcommand{\bi}{{\boldsymbol i}}
\newcommand{\bx}{{\boldsymbol p}}\newcommand{\bp}{{\boldsymbol p}}
\newcommand{\bpp}{{\boldsymbol P}}\newcommand{\bq}{{\boldsymbol q}}
\newcommand{\by}{{\boldsymbol q}}
\newcommand{\bz}{{\boldsymbol z}}
\newcommand{\eps}{{\varepsilon}}
\newcommand{\proofend}{\hfill$\Box$\bigskip}
\newcommand{\Int}{{\operatorname{Int}}}
\newcommand{\pr}{{\operatorname{pr}}}
\newcommand{\grad}{{\operatorname{grad}}}
\newcommand{\rk}{{\operatorname{rk}}}
\newcommand{\im}{{\operatorname{Im}}}
\newcommand{\sk}{{\operatorname{sk}}}
\newcommand{\const}{{\operatorname{const}}}
\newcommand{\Sing}{{\operatorname{Sing}}}
\newcommand{\conj}{{\operatorname{conj}}}
\newcommand{\Pic}{{\operatorname{Pic}}}
\newcommand{\Crit}{{\operatorname{Crit}}}
\newcommand{\Ch}{{\operatorname{Ch}}}
\newcommand{\discr}{{\operatorname{discr}}}
\newcommand{\Tor}{{\operatorname{Tor}}}
\newcommand{\Conj}{{\operatorname{Conj}}}
\newcommand{\val}{{\operatorname{val}}}
\newcommand{\Val}{{\operatorname{Val}}}
\newcommand{\res}{{\operatorname{res}}}
\newcommand{\add}{{\operatorname{add}}}
\newcommand{\tmu}{{\C\mu}}
\newcommand{\ov}{{\overline v}}\newcommand{\on}{{\overline n}}
\newcommand{\ox}{{\overline{x}}}
\newcommand{\tet}{{\theta}}
\newcommand{\Del}{{\Delta}}
\newcommand{\bet}{{\beta}}
\newcommand{\kap}{{\kappa}}
\newcommand{\del}{{\delta}}
\newcommand{\sig}{{\sigma}}
\newcommand{\alp}{{\alpha}}
\newcommand{\Sig}{{\Sigma}}
\newcommand{\Gam}{{\Gamma}}
\newcommand{\gam}{{\gamma}}
\newcommand{\Lam}{{\Lambda}}
\newcommand{\lam}{{\lambda}}
\newcommand{\SC}{{SC}}
\newcommand{\MC}{{MC}}
\newcommand{\nek}{{,...,}}
\newcommand{\cim}{{c_{\mbox{\rm im}}}}
\newcommand{\mathto}{\mathop{\to}}
\newcommand{\op}{{\overline p}}
\newcommand{\CL}{{\cal L}}
\newcommand{\CH}{{\cal CH}}
\newcommand{\w}{{\omega}}
\newcommand{\Sp}{{\mathbb S}}
\newcommand{\CX}{{\cal X}}
\newcommand{\LIC}{\text{\rm left}}
\newcommand{\CV}{{\cal V}}

\newcommand{\WW}{{W}}

\title{Welschinger invariants \\
of real Del Pezzo surfaces of degree $\ge 3$}
\author{Ilia Itenberg \and Viatcheslav Kharlamov \and Eugenii Shustin}

\date{}
\maketitle

\begin{abstract}
We give a recursive formula for purely real Welschinger invariants
of real Del Pezzo surfaces of degree $K^2\ge 3$, where
in the case
of surfaces of degree $3$
with two real components we introduce a certain modification
of Welschinger invariants and enumerate exclusively
the curves traced on the non-orientable component.
As an application, we prove the positivity of the
invariants under consideration and their logarithmic
asymptotic equivalence, as well as congruence modulo $4$,
to genus zero Gromov-Witten invariants.

\medskip\noindent {\bf MSC2010}:
Primary 14N10. Secondary 14P05, 14N35.

\medskip\noindent {\bf Keywords}: real rational curves,
enumerative geometry, Welschinger invariants, Caporaso-Harris
formula, cubic surfaces.
\end{abstract}


\bigskip
\bigskip

{\hskip3in From a dictionary for mathematicians:}

{\hskip3in {\bf Recursion} -- see {\it recursion}.}

\vskip10pt

{\sc\small \hfill "Mathematicians are also joking"}

\vskip-5pt

{\sc\small \hfill (compiled by S.~N.~Fedin)}

\section{Introduction}\label{intro}

In this paper we continue the study of purely real Welschinger
invariants of Del Pezzo surfaces. A particular interest of this
class of surfaces is related to the fact that the Welschinger
invariants of an unnodal real Del Pezzo surface are enumerative; in
particular, purely real Welschinger invariants of such a surface
count with certain signs the real rational curves that belong to a
given linear system and interpolate a suitable amount of real
points.

As we proved in \cite{IKS, IKS2,  IKS4, IKS6}, if $\Sigma$ is a real
Del Pezzo surface of degree $\geq 4$ with nonempty real part (except
the case of surfaces containing four disjoint $(-1)$-curves which
form two complex conjugate pairs) and $D$ is a nef and big real
divisor class, then the purely real Welschinger invariant
$W(\Sig,D)$ is positive (which implies the existence of
interpolating real rational curves). Furthermore, for these surfaces
the purely real Welschinger invariants and the corresponding genus
zero Gromov-Witten invariants are asymptotically equivalent in the
logarithmic scale, {\it i.e.},
\begin{equation}
\lim_{n\to+\infty}\frac{\log W(\Sig,nD)}{n\log
n}=\lim_{n\to+\infty} \frac{\log GW(\Sig,nD)}{n\log
n}=-DK_\Sig \label{eAS}
\end{equation}
(which implies a supexponential growth of the number
of interpolating real rational curves provided that
$\Sig$ is unnodal).

The main result of the present paper is a new recursive formula of
Caporaso-Harris type that applies to all purely real Welschinger
invariants of real Del Pezzo surfaces of degree $\ge3$ with nonempty
real part (Corollary \ref{c1}, section \ref{secRF}). Using this
formula, we extend the previous positivity and asymptotic results to
the plane blown up at $a$ real points and $b$ pairs of complex
conjugate points, where $a+2b\le6$, $b\le2$, 
as well as to the minimal
two-component real conic bundles over $\PP^1$ and the two-component real
cubic surfaces (see Theorems \ref{tn5}, \ref{t5}, and
Remark~\ref{rsss3}, section \ref{sec-pa}). Additionally, for the
surface~$\Sig$ which is the plane blown up at $a \leq 6$ real
points, we prove the monotone dependence of $W(\Sig, D)$ on the
divisor class $D$ and the Mikhalkin-type congruence
$$W(\Sig, D) = GW(\Sig, D)\mod4$$ (both claims were known
before for $a\le5$, {\it cf.} \cite{BM, IKS6}).

The present paper contains the first treatment of the positivity and
asymptotic relations of Welschinger invariants for surfaces having
at least two connected components of the real point set. The
original purely real Welschinger invariants are no more
unconditionally positive in such a case (see Remark~\ref{rsss3}). We
introduce some variations in the definition of Welschinger signs
that give us modified invariants (see details in section
\ref{W-invariants}), which are positive and do satisfy logarithmic
equivalence with genus zero Gromov-Witten invariants (Theorem
\ref{t5}, section \ref{sec-sss3}).

Unlike our previous works~\cite{IKS3, IKS6}, here we do not use
tropical geometry to derive the recursive formula. Instead, we
convert to a real form a complex Caporaso-Harris type formula
obtained in~\cite{MS} for the plane blown up at~$7$ points. 

A tropical calculation of purely real Welschinger invariants of the
plane blown up at $6$ real points was recently proposed by
E.~Brugall\'e~\cite{Bru}.

{\bf Acknowledgments}. We are grateful to the referee for careful
reading and helpful suggestions, and to E.~Brugall\'e for pointing
us a mistake in one of the table items of our Example \ref{ex}
and a missing deformation label in Section \ref{dedi}. 

A considerable part of the work on this text
was done during our visits to the {\it Max-Planck-Institut f\"{u}r
Mathematik} in Bonn in 2010, \emph{Mathematisches Forschungsinstitut
Oberwolfach} in 2011, and visits of the third author to {\it
Universit\'e de Strasbourg}. We thank these institutions for the
support and excellent working conditions.

The first two authors were partially funded by
the ANR-09-BLAN-0039-01 grant of {\it Agence Nationale de la
Recherche}, and are members of FRG Collaborative Research "Mirror
Symmetry \& Tropical Geometry" (Award No. $0854989$).
The third author enjoyed a support from the Israeli
Science Foundation grant no. 448/09 and from the
Hermann-Minkowski-Minerva Center for Geometry at the Tel Aviv
University.

\section{Welschinger invariants}\label{W-invariants}

Recall the original definition of Welschinger invariants in a form
adapted to the case of Del Pezzo surfaces. Let $\Sig$ be a real {\it
unnodal} ({\it i.e.}, not containing any rational $(-n)$-curve,
$n\ge 2$) Del Pezzo surface, and let $D\subset\Sigma$ be a real
effective divisor class. Consider a connected component $F$ of the
real point set $\R\Sigma$ of~$\Sig$ and a generic set~$\bp\subset F$
of $c_1(\Sig)\cdot D-1$ points. The set ${\cal R}(\Sig, D, \bp )$ of
real rational curves $C\in|D|$ passing through the points of~$\bp$
is finite, and all these curves are nodal and irreducible. Due to
the Welschinger theorem~\cite{W1} (and the genericity of the complex
structure on $\Sig$), the number
\begin{equation}
W(\Sig,D,\bp)=\sum_{C\in{\cal R}(\Sig,D, \bp)}(-1)^{s(C)}\ ,
\label{WS1}\end{equation} where $s(C)$ is the number of {\it
solitary nodes} of $C$ ({\it i.e.}, real points, where a local
equation of the curve can be written over $\R$ in the form
$x^2+y^2=0$), does not depend on the choice of a generic
set~$\bp\subset F$. We denote this (original) Welschinger invariant
by $\WW(\Sig,D,F)$.

If $\R\Sigma$ has more than one connected component (for example, if
$\Sigma$ is a cubic surface and $\R\Sigma$ has two connected
components), we modify the above construction of invariants in the
following way. Let, as above, $F$ be one of these components. For a
real nodal curve $C \subset \Sig$, we introduce its {\it Welschinger
weight reduced to} $F$ by putting $w_{F}(C)=(-1)^{s(C,F)}$, where
$s(C,F)$ is the number of real solitary nodes of $C$ belonging to
$F$. Then, given a real effective divisor class $D$ on $\Sig$, and a
generic set $\bz$ of $c_1(\Sig)D-1$ points in $F$, we define the
{\it Welschinger number of $D$ reduced to} $F$ by the formula
\begin{equation}W_{F}(\Sig,D,\bz)=\sum_{C\in{\cal R}(\Sig,D,
\bp)}w_{F}(C)\ .\label{WS2}\end{equation}

Such a {\it twisting} of the Welschinger construction can be
reformulated and slightly generalized. In addition to choosing one
of the real components, $F$, let us pick a homology class $\phi\in
H_2(\Sigma\setminus F; \Z/2\Z)$ invariant under the action of
complex conjugation, $\conj_*\phi=\phi$. Given a real effective
divisor class $D$ on $\Sig$, and a generic set $\bz$ of
$c_1(\Sig)D-1$ points in $F$, we define the {\it twisted Welschinger
number of $D$} by the formula
\begin{equation}W_{\phi}(\Sig,D,\bz)=\sum_{C\in{\cal R}(\Sig,D,
\bz)}w_\phi(C),\quad w_\phi(C)=(-1)^{s(C)+ C_\pm\circ\, \phi}\
,\label{ws}\end{equation} where $C_\pm$ denotes any of the two
halves of $C$ (any of the two discs cut from $C$ by $\R C$) and
$C_\pm\circ\, \phi\in\Z/2\Z$. Clearly, when $\phi_F$ is the homology
class realized in $H_2(\Sigma\setminus F; \Z/2\Z)$ by the union of
the components of $\R\Sigma\setminus F$, we get
$$
W_{\phi_F}(\Sig,D,\bz)=W_{F}(\Sig,D,\bz)\ .$$

\begin{proposition}\label{nl1}
The number $ W_{\phi}(\Sig,D,\bz)$ does not depend on the
choice of a generic set $\bz$ of $c_1(\Sig)D-1$ points in $F$.
\end{proposition}

{\bf Proof}. The statement is an immediate consequence of the
invariance of $\WW(\Sig,D,F)$ due to the following observation.

In a one-parametric family of curves $C(t)$
of class $D$ interpolating $c_1(\Sig)\circ D-2$ fixed generic points of $F$ and 
one additional point
of $F$ moving generically, the homology classes of the discs $C_\pm(t)$ are 
jumping only at those moments
$t=t_0$ when the curve $C(t_0)$ splits into two irreducible components $C'(t_0)$ 
and $C''(t_0)$. When such
a jump happens, the number $(-1)^{s(C(t))+ C_+(t)\circ\,\phi}$
does not change, since $C_+(t<t_0)\circ\,\phi=C'_+(t_0)\circ\,\phi+C''_+(t_0)\circ\,
\phi=C'_+(t_0)\circ\,\phi+\conj_*C''_+(t_0)\circ\,\conj_*\phi=
C'_+(t_0)\circ\,\phi+C''_-(t_0)\circ\,\phi= C_+(t>t_0)\circ\,\phi$.
\proofend

\smallskip

The above proposition implies existence of modified Welschinger invariants $\WW_\phi(\Sig,D)=W_\phi(\Sig,D,\bz)$..
As a particular example, we may take $\phi=\phi_F$, the fundamental class of the 
union of real components of
$\Sig$ different from $F$. This is the invariant which we use below
in the case of two-component real cubic surfaces
(and conic bundles);
we denote it, in accordance with our previous notation, $\WW_F(\Sig, D) $.

One may also choose as $\phi$ any combination of the fundamental
classes of real components different from $F$, and more generally,
combine them with vanishing classes between a pair of real
components different from $F$. In fact, one may prove, that for
multi-component real structures there do exist twists such that some
of the curves $C$ in ${\cal R}(\Sig,D, \bp)$ (with certain $D$
depending on the twist) change and some do not change the sign with
respect to the original Welschinger definition. Note also that the
number of independent possible twists is preserved under real
blow-ups (that is, a blow-up at a real point or at a pair of complex
conjugate points), so that interesting twists exist only for
surfaces with a disconnected real part.

\section{Recursive formula for Welschinger invariants}\label{Wel}

\subsection{Preliminaries}\label{sec1n}
In this chapter, we consider unnodal Del Pezzo surfaces of
degree~$3$. From the complex point of view, such a surface, denoted
by $\PP^2_6$, is the complex projective plane blown up at $6$ points
in general position. Denote by $L\subset\PP^2_6$ the strict
transform of a generic line, and by $E_1,...,E_6$ the exceptional
divisors of the blow up.

We equip $\Sig = \PP^2_6$ with a real structure, {\it i.e.}, an
anti-holomorphic involution $\conj:\Sig\to\Sig$. Then, as is well
known, the surface becomes isomorphic either to $\PP^2_{a,b}$,
$a+2b=6$, that is, the plane $\PP^2$ (equipped with its standard
real structure) blown up at $a$ real and $b$ pairs of complex
conjugate points, all in general position, or to $\CS$, a real cubic
surface with the two-component real part. By $F$ we denote the
non-orientable connected component of $\R\Sig$ (which is the only
one if $\Sig=\PP^2_{a,b}$), and we choose a class $\phi\in
H_2(\Sigma\setminus F; \Z/2\Z)$ invariant under the action of
complex conjugation, $\conj_*\phi=\phi$.

We pick a real smooth $(-1)$-curve $E$ on~$\Sig$ with $\R E\subset
F$: if $\Sig=\PP^2_{a,b}$, then choose $E\in|L-E_1-E_2|$,
where $E_1$ and $E_2$ are assumed to be either both real, or complex
conjugate;
if $\Sig = \CS$, then choose for $E$ any of the three lines whose real parts
are contained in $F$ (see, for example,~\cite{Seg}).

By $\Pic(\Sig,E)$ we denote the subsemigroup of $\Pic(\Sig)$
generated by (complex) irreducible curves, crossing $E$
non-negatively. The involution of complex conjugation~$\conj$ acts
on $\Pic(\Sig, E)$. By $\Pic^\R(\Sig,E)$ we denote the disjoint
union of the sets
$$
\left\{D\in\Pic(\Sig,E)\ :\ \conj D=D\right\}$$ and $$
\left\{\{D_1,D_2\}\in\Sym^2(\Pic(\Sig,E))\ :\ 
\conj D_1=D_2\right\}\ .
$$
For an element $\D\in\Pic^\R(\Sig,E)$, define $[\D]\in\Pic(\Sig,E)$
by
$$[\D]=\begin{cases}D,\quad & \text{if} \,\,
\D=D,\ \text{a divisor class},\\
D_1+D_2,\quad & \text{if} \,\, \D=\{D_1,D_2\},\ \text{a pair of divisor
classes}.\end{cases}$$

Let $\Z^\infty_+$ be the direct sum of countably many additive
semigroups \mbox{$\Z_+=\{m\in\Z\ |\ m\ge 0\}$} with the standard
basis
$$\theta_i=(\alp_1,\alp_2,...),\quad \alp_i=1,\quad \alp_j=0,\ j\ne
i\ .$$
For $\alp, \alp' \in \Z^\infty_+$,
the relation $\alp \geq \alp'$ means that
$\alp - \alp' \in\Z^\infty_+$.
For $\alp=(\alp_1,\alp_2,...)\in\Z^\infty_+$ put
$$
\|\alp\|=\sum_{k=1}^\infty\alp_k,\quad I\alp=\sum_{k=1}^\infty
k\alp_k,\quad I^\alp=\prod_{k=1}^\infty
k^{\alp_k},\quad\alp!=\prod_{k=1}^\infty\alp_k!\ ,\quad
$$ and for $\alp^{(0)},...,\alp^{(m)},\alp\in\Z_+^\infty$ such that
$\alp^{(0)}+...+\alp^{(m)}\le\alp$, put $$\left(\begin{matrix}\alp\\
\alp^{(0)},...,\alp^{(m)}\end{matrix}\right)
=\frac{\alp!}{\alp^{(0)}!...\alp^{(m)}!(\alp-\alp^{(0)}-...-\alp^{(m)})!}\
.$$
Introduce also the semigroup
$$\Z^{\infty,\odd}_+=\{\alp\in\Z^\infty_+\ :\ \alp_{2i}=0,\ i\ge
1\}\ .$$

For an element $\D\in\Pic^\R(\Sig,E)$ and a vector
$\bet\in\Z^\infty_+$, put
$$R_\Sig(\D,\bet)=-[\D](K_\Sig+E)+\|\bet\|-\begin{cases}1,\quad
& \text{if} \,\, \D=D,\ \text{a divisor class},\\
2,\quad & \text{if} \,\, \D=\{D_1,D_2\},\ \text{a pair of divisor
classes}.\end{cases}$$

\subsection{Families of real curves on $\Sig$}\label{Welsch-inv}
Let $\D\in\Pic^\R(\Sig,E)$, and vectors
$\alp,\bet^{\re},\bet^{\ima}\in\Z^\infty_+$ satisfy
\mbox{$I(\alp+\bet^{\re}+2\bet^{\ima})=[\D]E$}. Let
$\bp^\flat=\{p_{i,j}\ :\ i\ge1,\ 1\le j\le\alp_i\}$ be a sequence of
$\|\alp\|$ distinct real generic points on $E$. Such tuples
$(\D,\alp,\bet^{\re},\bet^{\ima},\bp^\flat)$ are called {\it
admissible}.

By
$V^\R_\Sig(\D,\alp,\bet^{\re},\bet^{\ima},\bp^\flat)$ we denote the
closure of the family of real reduced curves $C$
belonging to the linear system defined by $[\D]$ and such that
\begin{enumerate}
\item[(i)] if $\D=D$ is a divisor class, then $C\in|D|$ is an irreducible (over $\C$)
rational curve,
\item[(ii)] if $\D=\{D_1,D_2\}$ is a pair of divisors classes,

then $C=C_1\cup
C_2$, where $C_1\in|D_1|$, $C_2\in|D_2|$ are distinct, irreducible,
rational, conjugate imaginary curves;
\item[(iii)] $C\cap E$
consists of the points $\bp^\flat$ and $\|\bet^{\re}+2\bet^{\ima}\|$
other points; $\|\bet^{\re}\|$ of the latter points are real, and
the remaining points form $\|\bet^{\ima}\|$ pairs of complex
conjugate points;
\item[(iv)] $C$ has one local
branch at each of the points of $C\cap E$, and the intersection
multiplicities of~$C$ with $E$ are as follows:
\begin{itemize}
\item $(C\cdot E)_{p_{i,j}}=i$ for all $i\ge 1$, $1\le
j\le\alp_i$,
\item for each $i\ge 1$, there are precisely
$\bet^{\re}_i$ real points $q\in(C\cap E)\setminus\bp^\flat$ such
that $(C\cdot E)_q=i$;
\item for each $i\ge 1$, there are precisely
$\bet^{\ima}_i$ pairs $q,q'$ of complex conjugate points of $C \cap
E$ such that $(C\cdot E)_q = (C\cdot E)_{q'} = i$.
\end{itemize}
\end{enumerate}

\begin{lemma}\label{l1}
Let $(\D, \alp, \bet^\re, \bet^\ima, \bp^\flat )$ be an admissible
tuple.
If $V^\R_\Sig(\D,\alp,\bet^{\re},\bet^{\ima},\bp^\flat)$ is
nonempty, then $R_\Sig(\D,\bet^{\re}+2\bet^{\ima})\ge0$, and each
component of $V^\R_\Sig(\D,\bet^{\re},\bet^{\ima},\bp^\flat)$ has
dimension $\le R_\Sig(\D,\bet^{\re}+2\bet^{\ima})$. Moreover, a
generic element of any component of
$V^\R_\Sig(\D,\alp,\bet^{\re},\bet^{\ima},\bp^\flat)$ of dimension
$R_\Sig(\D,\bet^{\re}+2\bet^{\ima})$ is a nodal curve, nonsingular
along $E$.
\end{lemma}

{\bf Proof}. Follows from \cite[Propositions 2.1 and 2.2]{MS}. 
\proofend

\smallskip

Let $(\D, \alp, \bet^\re, \bet^\ima, \bp^\flat )$ be an admissible
tuple such that
$R_\Sig(\D,\bet^{\re}+2\bet^{\ima})\ge 0$. Pick a set $\bp^\sharp$
of $R_\Sig(\D,\bet^{\re}+2\bet^{\ima})$ generic points of
$F\setminus E$ and denote by
$V^\R_\Sig(\D,\alp,\bet^{\re},\bet^{\ima},\bp^\flat,\bp^\sharp)$ the
set of curves belonging to
$V^\R_\Sig(\D,\alp,\bet^{\re},\bet^{\ima},\bp^\flat)$
and passing through the points of
$\bp^\sharp$.

\begin{lemma}\label{l2}
Let $(\D, \alp, \beta^\re, \beta^\ima,
\bp^\flat)$
be an admissible tuple.

(1) if $\D=D$ is a divisor class, then
$V^\R_\Sig(D,\alp,\bet^{\re},\bet^{\ima},\bp^\flat,\bp^\sharp)$ is a
finite set of real nodal irreducible rational curves, nonsingular
along $E$;

(2) if $\D=\{D_1,D_2\}$ is a pair of divisor classes, then
$V^\R_\Sig(\D,\alp,\bet^{\re},\bet^{\ima},\bp^\flat,\bp^\sharp)$ is
nonempty only if $\alp=\bet^{\re}=0$, $R_\Sig(\D,2\bet^{\ima})=0$,
and $\bp^\flat=\bp^\sharp=\emptyset$; furthermore, in this case
the set $V^\R_\Sig(\D,\alp,\bet^{\re},\bet^{\ima},\bp^\flat,\bp^\sharp)$
is finite.
\end{lemma}

{\bf Proof}. By Lemma \ref{l1} we have to show only that
$R_\Sig(\D,2\bet^{\ima})=0$ is necessary for the nonemptyness of
$V^\R_\Sig(\D,0,0,\bet^{\ima},\emptyset,\bp^\sharp)$ with
$\D=\{D_1,D_2\}$. A curve $C\in
V^\R_\Sig(\D,0,0,\bet^{\ima},\emptyset,\bp^\sharp)$ splits in the
following way:
$$C=C_1\cup C_2,\quad C_1\in|D_1|,\quad C_2\in|D_2|,\quad \conj C_1=C_2\ ,$$ and,
by \cite[Proposition 2.1]{MS}, the component $C_1$ varies in a
family of complex dimension
$$-D_1(K_\Sig+E)+\|\bet^{\ima}\|-1=\frac{1}{2}R_\Sig(\D,2\bet^{\ima})\
.$$ Hence, a curve $C$ can match at most
$\frac{1}{2}R_\Sig(\D,2\bet^{\ima})$ generic points in $F\setminus
E$, and the claim follows. \proofend

\begin{lemma}[see, for example,~\cite{Seg}]\label{ln1}
(1) The linear system $|-(K_\Sig+E)|$ is of dimension~$1$ and
contains precisely two nonsingular curves $Q', Q''$ tangent to $E$,
and five reducible curves; each of the latter curves consists of two
distinct smooth $(-1)$-curves intersecting at one point.

(2) If $\Sig=\CS$, then
\begin{enumerate}
\item[(i)] one of the five reducible curves in the linear
system $|-(K_\Sig + E)|$ is formed by two real lines; each of the
other four reducible curves is formed by two complex conjugate lines
which intersect in one real point.
\item[(ii)] $\Sig$ has exactly three real lines,
and these lines generate the semigroup of real effective divisor
classes on $\Sig$.
\end{enumerate}\end{lemma}

The three real lines of $\Sig = \CS$ are denoted by
$L_1$, $L_2$, and $L_3$ (if the contrary is not explicitly stated,
we always assume that $E = L_1$).
The lines forming the four pairs of complex conjugate lines
are denoted by $L_j$, $j = 4$, $\ldots$, $11$, in a way
that, for any $i = 2, 3, 4, 5$, the lines
$L_{2i}$ and $L_{2i+1}$ are complex conjugate.

The next two lemmas follow from
\cite[Proposition 2.3]{MS} and Lemma \ref{ln1}.

\begin{lemma}\label{l3}
Let $\Sig=\PP^2_{a,b}$, $a+2b=6$, and $E\in|L-E_1-E_2|$. Then, among
the sets $V^\R_\Sig(\D,\alp,\bet^{\re},\bet^{\ima},\bp^\flat)$,
where $(\D, \alp, \beta^\re, \beta^\ima, \bp^\flat)$ is an
admissible tuple such that
$$
[\D]E > 0, \;\;\; R_\Sig(\D, \beta^\re + 2\beta^\ima) = 0,
$$
the only nonempty sets are as
follows:
\begin{enumerate}
\item[(1)]
in the case of a divisor class $\D=D$,
\begin{enumerate}\item[(1i)]
$V^\R_\Sig(E_i,0,\theta_1,0,\emptyset)$, where $i=1,2$, consists of
one element, provided that $E_1$ and $E_2$ are real;
\item[(1ii)]
$V^\R_\Sig(L-E_i-E_j,0,\theta_1,0,\emptyset)$, where $3\le i<j\le
6$, consists of one element, provided that $E_i$ and $E_j$ are
either both real or are complex conjugate;
\item[(1iii)]
$V^\R_\Sig(-(K_\Sig+E)-E_i,0,\theta_1,0,\emptyset)$, where $i=1,2$,
consists of one element, provided that $E_1$ and $E_2$ are real;
\item[(1iv)] $V^\R_\Sig(-(K_\Sig+E),0,\theta_2,0,\emptyset)$ consists
of two elements $Q',Q''$, provided that $Q'$ and $Q''$ are both
real;
\item[(1v)]
$V^\R_\Sig(-(K_\Sig+E),\theta_1,\theta_1,0,\bp^\flat)$ consists of
one element;
\item[(1vi)]
$V^\R_\Sig(-s(K_\Sig+E)+L-s_1E_1-s_2E_2-E_i,\alpha,0,0,\bp^\flat)$
consists of one element, if $s\ge 0$, $0\le s_1,s_2\le1$,
$s_1+s_2\le2s$, $3\le i \le6$, the divisor $E_i$ is real, and the
divisors $E_1,E_2$ are real whenever $s_1 \ne s_2$;
\item[(1vii)]
$V^\R_\Sig(-s(K_\Sig+E)-s_1E_1-s_2E_2+E_i,\alpha,0,0,\bp^\flat)$
consists of one element, if $s\ge 1$, $0\le s_1,s_2\le1$,
$s_1+s_2<2s$, $3\le i\le6$, the divisor $E_i$ is real, and the
divisors $E_1,E_2$ are real whenever $s_1 \ne s_2$;
\end{enumerate}
\item[(2)]
in the case of a pair $\D = \{D_1, D_2\}$ of divisor classes,
\begin{enumerate}
\item[(2i)] $V^\R_\Sig(\{E_1,E_2\},0,0,\theta_1,\emptyset)$ consists
of one element, provided that $E_1$ and $E_2$ are complex conjugate;
\item[(2ii)]
$V^\R_\Sig(\{L-E_i-E_j,L-E_i-E_k\},0,0,\theta_1,\emptyset)$, where
$\{i,j,k\}\subset\{3,4,5,6\}$, consists of one element, if $E_i$ is
real, and $E_j,E_k$ are complex conjugate;
\item[(2iii)]
$V^\R_\Sig(\{-(K_\Sig+E)-E_1,-(K_\Sig+E)-E_2\},0,0,\theta_1,\emptyset)$
consists of one element, provided that $E_1$ and $E_2$
are complex conjugate;
\item[(2iv)]
$V^\R_\Sig(\{L-E_i-E_j,L-E_k-E_l\},0,0,\theta_1,\emptyset)$, where
$\{i,j,k,l\}=\{3,4,5,6\}$, consists of one element, if $E_i,E_j$ and
$E_k,E_l$ are two pairs of complex conjugate exceptional divisors;
\item[(2v)]
$V^\R_\Sig(\{-(K_\Sig+E),-(K_\Sig+E)\},0,0,\theta_2,\emptyset)$
consists of one element $\{Q',Q''\}$, provided that $Q'$ and $Q''$
are complex
conjugate.
\end{enumerate}\end{enumerate}
\end{lemma}

\begin{lemma}\label{l4}
Let $\Sig=\CS$ and $E=L_1$. Then, among the sets
$V^\R_\Sig(\D,\alp,\bet^{\re},\bet^{\ima},\bp^\flat)$, where $(\D,
\alp, \beta^\re, \beta^\ima, \bp^\flat)$ is an admissible tuple such
that
$$
[\D]E > 0, \;\;\; R_\Sig(\D, \beta^\re + 2\beta^\ima) = 0,
$$
the only nonempty sets are as
follows:
\begin{enumerate}\item[(1)] in the case of
a divisor class $\D=D$,
\begin{enumerate}\item[(1i)]
$V^\R_\Sig(L_i,0,\theta_1,0,\emptyset)$, where $i=2,3$, consists of
one element $L_i$; \item[(1ii)]
$V^\R_\Sig(-(K_\Sig+E),0,\theta_2,0,\emptyset)$ consists of two
elements $Q'$ and $Q''$;
\item[(1iii)] $V^\R_\Sig(-(K_\Sig+E),\theta_1,\theta_1,0,\bp^\flat)$,
consists of one element;

\end{enumerate}
\item[(2)] in the case of a pair $\D = \{D_1, D_2\}$ of
divisor classes,
\begin{enumerate}\item[]
$V^\R_\Sig(\{L_{2i},L_{2i+1}\},0,0,\theta_1,\emptyset)$, where
$i=2,3,4,5$, consists of one element $\{L_{2i},L_{2i+1}\}$.
\proofend
\end{enumerate}
\end{enumerate}
\end{lemma}

\subsection{Deformation diagrams}\label{dedi}

Let $(D,\alp,\bet^{\re},\bet^{\ima},\bp^\flat)$ be an admissible
tuple, where $D\in\Pic^\R(\Sig,E)$ is a divisor class and
$R_\Sig(D,\bet^{\re}+2\bet^{\ima})>0$. Pick a set
$\widetilde\bp^\sharp$ of $R_\Sig(D,\bet^{\re}+2\bet^{\ima})-1$
generic real points of $F\setminus E$, a generic real point $p\in
E\setminus\bp^\flat$, and a smooth real algebraic curve germ $\Lam$
crossing $E$ transversally at $p$. Denote by $\Lam^+=\{p(t)\ :\
t\in(0,\eps)\}$ a parameterized connected component of
$\Lam\setminus\{p\}$ with $\lim_{t\to0}p(t)=p$. There exists
$\eps_0>0$ such that, for all $t\in(0,\eps_0]$, the sets
$V_\Sig(D,\alp,\bet^{\re},\bet^{\ima},\bp^\flat,
\widetilde\bp^\sharp\cup\{p(t)\})$ are finite, their elements remain
nodal, nonsingular along $E$ as $t$ runs over the interval
$(0,\eps_0]$, and the closure in
$V_\Sig(D,\alp,\bet^{\re},\bet^{\ima},\bp^\flat)$ of the family
\begin{equation}V=\bigcup_{t\in(0,\eps_0]}V_\Sig(D,\alp,
\bet^{\re},\bet^{\ima},\bp^\flat,
\widetilde\bp^\sharp\cup\{p(t)\})\label{eDD}\end{equation} is a
union of real algebraic arcs which are disjoint for $t>0$. This
closure $\overline V$ is called a \emph{deformation diagram} of
$(D,\alp,\bet^{\re},\bet^{\ima},\bp^\flat,\widetilde\bp^\sharp,p)$.
The elements of
$V_\Sig(D,\alp,\bet^\re,\bet^\ima,\bp^\flat,\widetilde\bp^\sharp\cup
\{p(\eps_0)\})$ are called \emph{leaves} of the deformation diagram,
and the elements of $\overline V\setminus V$ are called \emph{roots}
of the deformation diagram.

\begin{lemma}\label{DC}
Each connected component of a deformation diagram $\overline V$
defined by (\ref{eDD}) contains exactly one root. The roots are
curves $C\in|D|$ of the following two types.

(I) The curve $C$ is a generic member of an
$(R_\Sig(D,\bet^{\re}+2\bet^{\ima})-1)$-dimensional component of one
of the families
$$V^\R_\Sig(D,\alp+\theta_j,\bet^{\re}-\theta_j,\bet^{\ima},
\bp^\flat\cup\{p\}, \widetilde\bp^\sharp\})\ ,$$ where $j$ is a
natural number such that $\bet^{\re}_j>0$.

(R) The curve $C$ decomposes into $E$ and curves
of the following four types (for each type, the collection
of curves can be empty):
\begin{enumerate}
\item[(R1)] distinct reduced irreducible over $\R$ curves $C^{(i)}$,
$1\le i\le m$, which are generic members in some
$R_\Sig(\D^{(i)},(\bet^{\re})^{(i)}+2(\bet^{\ima})^{(i)})$-dimensional
components of families
$$V^\R_\Sig(\D^{(i)},\alp^{(i)},(\bet^{\re})^{(i)},(\bet^{\ima})^{(i)},
(\bp^\flat)^{(i)},(\bp^\sharp)^{(i)})\ ,$$
respectively, where $\D^{(i)}$ is a divisor class if $C^{(i)}$ is
irreducible over $\C$, and is a pair of divisor classes if $C^{(i)}$
is the union of two complex conjugate components, and, in
addition, $\D^{(i)}$ is neither $-(K_\Sig+E)$, nor
$\{-(K_\Sig + E), -(K_\Sig + E)\}$,
$1\le i\le m$,
\item[(R2)] distinct
curves $jQ(p_{j,s})$, where $p_{j,s}$ runs over some subset
$(\bp^\flat)^{(0)}$ of $\bp^\flat$, and $Q(p_{j,s})\in|-(K_\Sig+E)|$
is the {\rm (}real\/{\rm )} curve passing through $p_{j,s}$,
\item[(R3)] distinct curves $j(z)Q(z)$, where $z$ runs
over some subset $(\bp^\sharp)^{(0)}$ of $\widetilde\bp^\sharp$,
$j(z)\ge1$, and $Q(z)\in|-(K_\Sig+E)|$ is the {\rm (}real\/{\rm )}
curve passing through $z$,
\item[(R4)] curves $l'Q'$ and $l''Q''$, where $Q',Q''\in|-(K_\Sig+E)|$
are the two curves tangent to $E$ ({\it cf.} Lemma \ref{ln1}(1)),
and $l'=l''\ge0$ if $Q',Q''$ are complex conjugate, and $l',l''\ge0$
if $Q',Q''$ are real.
\end{enumerate}
Furthermore, the parameters of the above decomposition
are subject to the following restrictions:
\begin{itemize}
\item
$\sum_{i=1}^m\alp^{(i)}\le\alp-\alp^{(0)}$, where $\alp^{(0)}$
encodes the sequence of multiplicities $j$ over all
$p_{j,k}\in(\bp^\flat)^{(0)}$,
\item $\sum_{i\in S}(\bet^{\ima})^{(i)}=\bet^{\ima}$,
where $S=\{i\in[1,m]\ :\ \D^{(i)}\ \text{is a divisor class}\}$,
\item
$\bet^{(0)}\le\bet^\re$, where $\bet^{(0)}$ encodes the sequence of
multiplicities $j(z)$ over all $z\in(\bp^\sharp)^{(0)}$,
\item there
is a sequence of vectors $\gam^{(i)}\in\Z_+^\infty$, $i\in S$, such that
$\|\gam^{(i)}\|=1$, $\gam^{(i)}\le(\bet^\re)^{(i)}$, $i\in S$, and
$\sum_{i\in S}((\bet^\re)^{(i)}-\gam^{(i)})=\bet^\re-\bet^{(0)}$,
\item
$\sum_{i=1}^m[\D^{(i)}]=D-E+(I\alp^{(0)}+I\bet^{(0)}+l'+l'')(K_\Sig+E)$.
\end{itemize}
\end{lemma}

{\bf Proof}. All claims follow from \cite[Proposition 2.6]{MS}. We
only make two comments. In the case (I), an imaginary moving
intersection point with $E$ cannot merge to $p$, since otherwise the
conjugate moving intersection point must merge to $p$ too. In the
case (R), in view of \cite[Lemma 2.9]{MS} and due to the
rationality of curves in
$V^\R_\Sig(D,\alp,\bet^{\re},\bet^{\ima},\bp^\flat,
\widetilde\bp^\sharp \cup
\{p(\eps_0)\}$,
{\it cf.} \cite[Corollary 2.20]{MS},
in the deformation of $C$ induced by $\overline V$, for
each irreducible over $\C$ component $C'\ne E$ of $C$, precisely one
of the intersection points of~$C'$ with $E\setminus\bp^\flat$ is
smoothed out, whereas the remaining intersection points of~$C'$ with
$E$ turn into smooth points of the deformed curve, where it crosses
$E$ with the same multiplicity as $C'$ (these intersection points
stay fixed if they were in $\bp^\flat$ or move along $E$ otherwise).
To get restrictions on the parameters of the decomposition, we
notice also that, by Lemmas \ref{l2}(2), \ref{l3}(2), and
\ref{l4}(2), if $C^{(i)}$ has a pair of complex conjugate
components, then each component crosses $E$ at a unique point, and
this point is imaginary. \proofend

\begin{lemma}\label{ldediI}
Let~$C$ be the root of a connected component $\delta$ of the
deformation diagram $\overline V$. Assume that $C \in V^\R_\Sig(D,
\alp + \theta_j, \beta^\re - \theta_j, \beta^\ima, \bp^\flat \cup
\{p\}, \widetilde\bp^\sharp)$ is of type (I). If $j$ is odd, then
$\delta$ has a unique leaf; if $j$ is even, then $\delta$ has two
leaves. In both cases, each leaf has the same real topology as the
root $C$.
\end{lemma}

{\bf Proof}. Straightforward from \cite[Formulas (23) and (24)]{MS}.
\proofend

Let $C$ be the root of a connected component $\delta$ of the
deformation diagram $\overline V$. Assume that~$C$  is of type (R).
To describe the leaves of~$\delta$, we introduce \emph{deformation
labels}, certain nodal curves specified below. Each deformation
label is given by a polynomial equation $\psi(x,y)=0$ in the toric
surface $\Tor(\psi)$ defined by the Newton polygon of $\psi$. The
following list contains the deformation labels and, in the case of
real deformation labels, specifies the number of their solitary
nodes.
\begin{enumerate}
\item[(DL1)]
Two deformation
labels defined by the equations
$$\psi_1(x,y)=y^2+1+y\cdot\cheb_j(x)=0\quad\text{and}
\quad\psi_2(x,y)=\psi_1(\sqrt{-1} \ x,\sqrt{-1}y)=0\ ,$$
where $j$ is a positive odd number and $\cheb_j(\cos x)=\cos jx$ are the 
Chebyshev polynomials; the former curve
has $j - 1$ solitary nodes, the latter curve has no
real nodes. 
\item[(DL2)]
Two deformation
labels defined by the equations
$$\psi_1(x,y)=y^2+1+y\cdot\cheb_j(x)=0\quad\text{and}
\quad\psi_2(x,y)=\psi_1(\sqrt{-1} \ x,y)=0\ ,$$
where~$j$ is a positive even number; the former curve
has $j - 1$ solitary nodes, the latter curve has no
solitary nodes.
\item[(DL3)]
Deformation label defined by the equation
$$\psi(x,y)=(x-1)(y^j-x)=0\ ;$$
this curve has no solitary nodes.
\item[(DL4)]
Deformation label defined by the equation
$$\psi(x,y)=(x-1)(1+x((y+2^{1/j})^j-1))=0\ ,$$
where $j$ is a positive odd number; this curve has no solitary nodes.
\item[(DL5)]
Two deformation labels defined by the equations
$$\psi(x,y)=(x-1)(1+x((y\pm 2^{1/j})^j-1))=0\ ,$$
where $j$ is a positive even number; these curves have no solitary nodes.
\item[(DL6)]
Deformation label defined by the equation
$$\psi(x,y)=1+\frac{y+x^2}{2y}\left(\cheb_{l'+1}
\left(\frac{\pm y}{2^{(l'-1)/(l'+1)}}+y'\right)-1\right)=0\ ,$$
where $l'$ is a positive even number, and $y'$ is the only positive
simple root of $\cheb_{l' + 1}(y) - 1$; this curve has either $l'$
solitary nodes, or no solitary nodes at all.
\item[(DL7)]
Two deformation labels given by the equations
$$\displaylines{
\psi_1(x,y)=1+\frac{y+x^2}{2y}\left(\cheb_{l'+1}
\left(\frac{y}{2^{(l'-1)/(l'+1)}} + y'\right)-1\right)=0\ , \cr
\psi_2(x,y) = 1 + \frac{y + x^2}{2y}\left(\cheb_{l' + 1}
\left(\frac{y}{2^{(l' - 1)/(l' + 1)}} - y'\right) - 1\right) = 0,
}
$$
where $l'$ is a positive odd number, and $y'$ is the only positive
simple root of $\cheb_{l'+1}(y)-1$; the former curve has $l'$
solitary nodes, the latter curve has no solitary nodes.
\item[(DL8)]
$l'+1$ deformation labels defined by the equations
$$\psi(x,y)=1+\frac{y+\sqrt{-1} \ x^2}{2y}\left(\cheb_{l'+1}
\left(\frac{y\eps}{2^{(l'-1)/(l'+1)}}+y'\right)-1\right)=0,\quad
\eps^{l'+1}=1\ ,$$ where $l'$ is a positive integer.
\item[(DL9)]
$l'' + 1$ deformation labels defined by the equations
$$\psi(x,y)=1+\frac{y - \sqrt{-1} \ x^2}{2y}\left(\cheb_{l''+1}
\left(\frac{y\eps}{2^{(l''-1)/(l''+1)}}+y'\right)-1\right)=0,\quad
\eps^{l''+1}=1\ ,$$
where $l''$ is a positive integer.
\end{enumerate}

Consider now
the
following data:
\begin{enumerate}
\item[(C1)] choose
a sequence of vectors
$\gam^{(i)}\in\Z_+^\infty$, $i = 1$, $\ldots$, $m$, such that
\begin{itemize}
\item
$\|\gam^{(i)}\|=1$ for $i\in S=\{i\ :\ 1\le i\le m,\ \D^{(i)}\
\text{is a divisor class}\}$,
\item $\gam^{(i)}=0$ for $1 \leq i \leq m$ and $i\not\in S$,
\item
$\gam^{(i)}\le(\bet^\re)^{(i)}$, $i = 1$, $\ldots$, $m$,
\item
$\sum_{i = 1}^m((\bet^\re)^{(i)}-\gam^{(i)})=\bet^\re-\bet^{(0)}$;
\end{itemize}
\item[(C2)]
for each $i \in S$, choose a real point $q_i\in (C^{(i)}\cap
E)\setminus\bp^\flat$ such that $(C^{(i)}\cdot E)_{q_i}=j$, where
$\gam^{(i)}=\theta_j$;
\item[(C3)]
for each $z \in (\bp^\sharp)^{(0)}$, choose a point $q(z)$ which is
of one of the two real points of $Q(z)\cap E$;
\item[(C4)]
now to each chosen in (C2) real point $q_i$, $i\in S$, we assign
a number $\eps_i=\pm1$ in the following way: take local conjugation-invariant coordinates
$x,y$ in a neighborhood of $q_i$ so that $q_i=(0,0)$, $E=\{y=0\}$, and
$C=y(y+x^j+\text{h.o.t.})$; then (see \cite[Formulas (25), (28), and (30)]{MS})
any curve in the family $V$ defined in
(\ref{eDD}) is given by $y(y+x^j+\text{h.o.t.}+O(t^{>0}))+t(a+O(t,x))=0$, where
$t>0$ is a regular parameter on $\Lambda^+$,
$a\ne0$, and $\eps_i=\sign(a)$ depends only on the relative position of the points in $\widetilde\bp^\sharp\cup\{p\}$
and the points chosen in (C2), (C3) on $\R E$.
\end{enumerate}

Denote by ${\mathcal C}(C)$ the set of all possible choices of data
(C1)-(C4).
For any element of ${\mathcal C}(C)$, we construct a {\it suitable
deformation label collection}, which is a sequence of deformation labels.
If there is a point $q_i$ chosen in (C2) with an even
$j=
(C^{(i)}\cdot E)_{q_i}$
and $\eps_i=-1$, we define a deformation label collection to
be empty, otherwise, we construct it as follows
\begin{itemize}
\item 
one deformation label of type (DL1) for each point
$q_i$
chosen in (C2) with an odd $j=(C^{(i)}\cdot E)_{q_i}$; if $\eps_i=1$, we take $\psi_1(x,y)=0$, if
$\eps_i=-1$, we take $\psi_2(x,y)=0$;
\item any of the two deformation labels of type (DL2) for each point
$q_i$ chosen in (C2) with an even $j=(C^{(i)}\cdot E)_{q_i}$ and $\eps_i=1$;
\item one deformation label of type (DL3)
for each component $jQ(p_{j,s})$ of $C$,
\item one deformation label for each
component $j(z)Q(z)$ of $C$ and each point $q(z)$
chosen in (C3); this deformation label is of type (DL4) or (DL5)
depending on the parity of $j = j(z)$;
\item if $Q'$ and $Q''$ are real,
one deformation label for $l'Q'$ and
one deformation label for $l''Q''$,
the former (respectively, the latter)
deformation label is of type (DL6) or (DL7)
depending on the parity of $l'$ (respectively, $l''$);
\item if $Q'$ and $Q''$ are complex conjugate
(in this case $l' = l''$), one deformation label of type (DL8) for
$l'Q'$ and the complex conjugate deformation label of type (DL9) for
$l''Q''$.
\end{itemize}
Denote by ${\mathcal Def}(C,\sigma)$ the set of all suitable
deformation label collections of a given element $\sigma\in{\mathcal
C}$.

\begin{lemma}\label{lDP}
Let $C$ be the root of a connected component $\delta$ of the
deformation diagram $\overline V$. Assume that~$C$ is of type (R).

(1) Suppose that $\bet^{(0)}\in\Z_+^{\infty,\odd}$,
$\sum_{i=1}^m(\bet^{\re})^{(i)}-\bet^\re\in\Z_+^{\infty,\odd}$, and
either $Q',Q''$ are complex conjugate, or $Q',Q''$ are real and
$l',l''$ are both even. Then, there is a one-to-one correspondence
between the set of leaves of $\delta$ and the disjoint union of the
sets ${\mathcal Def}(C,\sigma)$ over all $\sigma\in{\mathcal C}(C)$.

Suppose that either $\bet^{(0)}\not\in\Z_+^{\infty,\odd}$, or
$\sum_{i=1}^m(\bet^{\re})^{(i)}-\bet^\re\not\in\Z_+^{\infty,\odd}$,
or both $Q',Q''$ are real and at least one of $l',l''$ is odd. Then,
the set of leaves of $\delta$ is in one-to-one correspondence with
the disjoint union of sets ${\mathcal Def}(C,\sigma)$, where
$\sigma$ runs over some nonempty subset of ${\mathcal C}(C)$.

(2) The set of solitary nodes of each leaf of $\delta$ bijectively
corresponds to the disjoint union of the sets of solitary nodes of
the corresponding deformation labels and the sets of solitary nodes
of the components $C^{(i)}$, $i = 1$, $\ldots$, $m$, of $C$.
Moreover, the solitary nodes coming from the deformation labels all
belong to the connected component $F \subset \R\Sig$ which contains
the line $\R E$.
\end{lemma}

{\bf Proof}. Statement (1) follows from \cite[Lemma 2.19]{MS}
(one-to-one correspondence) and \cite[Lemma 2.9]{MS} (geometry of
deformation), which describe all complex deformations of $C$ via
so-called deformation patterns. Our deformation labels can be viewed
as normalized versions of these deformation patterns. The restricted
correspondence in the second case of assertion (1) comes from the
fact that, for some elements of ${\mathcal C}(C)$, all complex
deformations of $C$ are non-real. Statement (2) follows from
\cite[Lemmas 2.10, 2.12, 2.14, and 2.15]{MS}, where one can find a
complete description of complex deformation patterns and formulas
for them. \proofend

\subsection{Welschinger numbers}\label{sec5}
For any admissible tuple $(\D, \alp, \bet^\re, \bet^\ima, \bp^\flat
)$ such that $R_\Sig(\D,\bet^{\re}+2\bet^{\ima})\ge 0$, and for any
set $\bp^\sharp$ of $R_\Sig(\D,\bet^{\re}+2\bet^{\ima})$ generic
points of $F\setminus E$, consider the set
$V^\R_\Sig(\D,\alp,\bet^{\re},\bet^{\ima},\bp^\flat,\bp^\sharp)$,
see section \ref{Welsch-inv}. According to Lemma \ref{l2}, this set
is finite and consists of real nodal irreducible rational curves.
Put
\begin{equation}
W_{\Sig,\phi}(\D,\alp,\bet^{\re},\bet^{\ima},
\bp^\flat,\bp^\sharp)=\sum_{C\in
V^\R_\Sig(\D,\alp,\bet^{\re},\bet^{\ima},\bp^\flat,\bp^\sharp)}(-1)^{s(C)+
C_\pm\circ\, \phi}\ .\label{e1}
\end{equation}

In view of Proposition
\ref{nl1}, for any divisor class $D\in\Pic^\R(\Sig,E)$ and a set
$\bp^\sharp$ of $c_1(\Sig)D-1$ distinct generic points of
$F\setminus E$, one has
\begin{equation}\WW_\phi(\Sig,D)=\sum_{k+2l=DE}W_{\Sig,\phi}
(D,0,k\theta_1,l\theta_1,\emptyset,\bp^\sharp)\ .\label{e35}\end{equation}

Pick a divisor class $D_0\in\Pic^\R(\Sig,E)$, and put $N=\dim|D_0|$.
Note that the set
$$\Prec(D_0)=\{D\in\Pic^\R(\Sig,E)\ :\ D_0\ge D\}$$ is finite, and we have $\dim|D|\le N$
for each $D\in\Prec(D_0)$. Furthermore, for each nonempty variety
$V^\R_\Sig(D,\alp,\bet^{\re},\bet^{\ima},\bp^\flat)$ with
$D\in\Prec(D_0)$, we have
$$\|\alp\|+R_\Sig(D,\bet^{\re}+2\bet^{\ima})\le N\ .$$

\begin{lemma}\label{D0}
Let $D_0\in \Pic^\R(\Sig,E)$ be a divisor class with
$N=\dim|D_0|>0$. Then, there exists a sequence
$\Lam(D_0)=(\Lam_i)_{i=1,...,N}$ of $N$ disjoint smooth real
algebraic arcs in $\Sig$, which are parameterized by
$t\in[0,1]\mapsto p_i(t)\in\Lam_i$, such that $p_i(0)\in E$,
$i = 1$, $\ldots$, $N$,
the
arcs $\Lam_i$ are transverse to $E$ at $p_i(0)$,
$i = 1$, \ldots, $N$, and
the following condition holds:

for an arbitrary admissible tuple
$(D,\alp,\bet^{\re},\bet^{\ima},\bp^\flat)$, disjoint subsets
$J^\flat,J^\sharp\subset\{1,...,N\}$, and a positive integer $k\le
N$ such that
\begin{enumerate}\item[(i)] $D\le D_0$, \item[(ii)] $R_\Sig(D,\bet^{\re}+2\bet^{\ima})>0$,
\item[(iii)] $i<k<j$ for all $i\in J^\flat$, $j\in J^\sharp$, \item[(iv)] the
number of elements in $J^\sharp$ is equal to
$R_\Sig(D,\bet^{\re}+2\bet^{\ima})-1$,
\item[(v)] $\bp^\flat=\{p_i(0)\ :\ i\in J^\flat\}$, \end{enumerate} the
closure of the family
$$\bigcup_{t\in(0,1]}V_\Sig(D,\alp,\bet^{\re},\bet^{\ima},\bp^\flat,
\widetilde\bp^\sharp\cup\{p_k(t)\})\ ,$$ where $\widetilde\bp^\sharp=
\{p_j(1)\}_{j\in J^\sharp}$, is a deformation diagram of
$(D,\alp,\bet^{\re},\bet^{\ima},\bp^\flat,\widetilde\bp^\sharp,p_k(0))$.
\end{lemma}

{\bf Proof}. Take a sequence $\widehat\Lam_i$, $i=1,...,N$, of
disjoint smooth real algebraic arcs in $\Sig$, which are
parameterized by $t\in[0,1]\mapsto p_i(t) \in \widehat\Lam_i$, such that
$(p_i(0))_{i=1,...,N}$ is a generic sequence of points in $E$,
and the arcs $\widehat\Lam_i$ are transverse to $E$ at $p_i(0)$,
$i=1,...,N$. We will inductively shorten and reparameterize these
arcs in order to satisfy the diagrammatic condition.

Suppose that we have constructed the arcs $\Lam_1,...,\Lam_{k-1}$,
$1\le k\le N$. There are finitely many admissible tuples
$(D,\alp,\bet^{\re},\bet^{\ima},\bp^\flat)$ and subsets
$J^\flat,J^\sharp\subset\{1,...,N\}$ satisfying restrictions (i)-(v)
in Lemma. Given such data
$D,\alp,\bet^{\re},\bet^{\ima},\bp^\flat,J^\flat,J^\sharp$, we take
a closed neighborhood
$\Lam_k(D,\alp,\bet^{\re},\bet^{\ima},\bp^\flat,J^\flat,J^\sharp)$
of $p_k(0)$ in $\widehat\Lam_k$, parameterized by $[0,\eps_k]$, such
that the closure of the family
$$\bigcup_{p'\in\Lam_k(D,\alp,\bet^{\re},\bet^{\ima},\bp^\flat,J^\flat,J^\sharp),p'\ne p_k(0)}
V_\Sig(D,\alp,\bet^{\re},\bet^{\ima},\bp^\flat,\widetilde\bp^\sharp\cup\{p'\})\
,$$ where $\widetilde\bp^\sharp=\{p_i(\eps_i)\}_{1\le i<k}$, is a
deformation diagram of
$(D,\alp,\bet^{\re},\bet^{\ima},\bp^\flat,\widetilde\bp^\sharp,p_k(0))$.
Then we define
$$\Lam_k=\bigcap_{(D,\alp,\bet^{\re},\bet^{\ima},\bp^\flat,J^\flat,J^\sharp)}
\Lam_k(D,\alp,\bet^{\re},\bet^{\ima},\bp^\flat,
J^\flat,J^\sharp)$$ and reparameterize this arc by $[0,1]$.
\proofend

Take a divisor class $D_0\in\Pic^\R(\Sig,E)$ such that
$N=\dim|D_0|>0$ and a sequence of arcs $(\Lam_i)_{i=1,...,N}$ as in
Lemma \ref{D0}. Given two subsets
$J^\flat,J^\sharp\subset\{1,...,N\}$ such that $i<j$ for all $i\in
J^\flat$, $j\in J^\sharp$, we say that the pair of point
configurations
$$\bp^\flat=\{p_i(0)\ :\ i\in J^\flat\},\quad
\bp^\sharp=\{p_j(1)\ :\ j\in J^\sharp\}$$ is \emph{in $D_0$-CH
position}.

\begin{proposition}\label{l12}
Fix a tuple $(D, \alp, \bet^\re, \bet^\ima)$, where~$D \in
\Pic^\R(\Sig, E)$ is a divisor class,
$\alp,\bet^\re\in\Z_+^{\infty,\odd}$, and $\bet^\ima \in
\Z_+^{\infty}$ such that $R_\Sig(D, \bet^\re + 2\bet^\ima) > 0$.
Choose two point
sequences $\bp^\flat$ and $\bp^\sharp$ satisfying the following
restrictions:
\begin{itemize}
\item the tuple $(D, \alp, \bet^\re, \bet^\ima, \bp^\flat)$
is admissible,
\item the number of points in $\bp^\sharp$ is equal to $R_\Sig(D, \bet^\re + 2\bet^\ima)$,
\item the pair $(\bp^\flat, \bp^\sharp)$ is in $D_0$-CH position
for some divisor class $D_0 \in \Pic^\R(\Sig, E)$, $D_0 \geq D$.
\end{itemize}
Then, the number $W_{\Sig,\phi}(D, \alp, \bet^\re, \bet^\ima,
\bp^\flat, \bp^\sharp)$ does not depend on the choice of sequences
$\bp^\flat$ and $\bp^\sharp$.
\end{proposition}

The proof is presented in section \ref{secn6}.

\smallskip

We skip $\bp^\flat$ and $\bp^\sharp$ in the notation of the above
numbers and simply write
$W_{\Sig,\phi}(D,\alp,\bet^{\re},\bet^{\ima})$ calling them {\bf
Welschinger numbers}.

\begin{proposition}\label{ini1}
Let $\Sig=\PP^2_{a,b}$, $a+2b=6$, and $E\in|L-E_1-E_2|$. If $\phi =
0 \in H_2(\Sigma\setminus F; \Z/2\Z)$, then among the Welschinger
numbers $W_{\Sig,
\phi}(\D,\alp,\bet^{\re},\bet^{\ima},\bp^\flat,\emptyset)$, where
$(\D,\alp,\bet^\re,\bet^\ima,\bp^\flat)$ is an admissible tuple such
that
$$
\alpha,\beta^{\re}\in\Z_+^{\infty,\odd}, \;\;\; [\D]E > 0,  \,\text{and} \;\;\;
R_\Sig(\D, \beta^\re + 2\beta^\ima) = 0,\,
$$
the only non-zero numbers
are as follows:
\begin{enumerate}
\item[(1)]
in the case of a divisor class $\D = D$,
\begin{enumerate}
\item[(1i)]
$W_{\Sig, \phi}(E_i,0,\theta_1,0,\emptyset,\emptyset) = 1$, where
$i=1,2$, and $E_1,E_2$ are real;
\item[(1ii)]
$W_{\Sig, \phi}(L-E_i-E_j,0,\theta_1,0,\emptyset,\emptyset) = 1$,
where $3\le i<j\le 6$, and $E_i,E_j$ are both real or are complex
conjugate;
\item[(1iii)]
$W_{\Sig, \phi}(-(K_\Sig+E)-E_i,0,\theta_1,0,\emptyset,\emptyset) =
1$, where $i=1,2$, and $E_1,E_2$ are real;
\item[(1iv)]
$W_{\Sig, \phi}(-(K_\Sig+E),\theta_1,\theta_1,0,\bp^\flat,\emptyset) = 1$;
\item[(1v)]
$W_{\Sig, \phi}(-s(K_\Sig + E) + L - s_1E_1 - s_2E_2 - E_i, \alp, 0,
0,\bp^\flat,\emptyset) = 1$, where $s \ge 0$, $0\le s_1,s_2\le 1$,
$s_1+s_2\le2s$, $3 \leq i \leq 6$, the sequence $\alp \in
\Z^{\infty, \odd}$ verifies $I\alp = 2s+1-s_1-s_2$, the divisor
$E_i$ is real, and the divisors~$E_1$ and~$E_2$ are both real if
$s_1 \ne s_2$;
\item[(1vi)]
$W_{\Sig, \phi}(-s(K_\Sig + E) - s_1E_1 - s_2E_2 + E_i, \alp, 0,
0,\bp^\flat,\emptyset) = 1$, where $s\ge 1$, $0\le s_1,s_2 \le1$,
$s_1+s_2<2s$, $3 \leq i \leq 6$, the sequence $\alp \in \Z^{\infty,
\odd}$ verifies $I\alp=2s-s_1-s_2$, the divisor $E_i$ is real, and
the divisors~$E_1$ and~$E_2$ are both real if $s_1 \ne s_2$;
\end{enumerate}
\item[(2)]
in the case of a pair $\D = (D_1, D_2)$ of divisor classes,
\begin{enumerate}
\item[(2i)] $W_{\Sig, \phi}(\{E_1,E_2\},0,0,\theta_1,
\emptyset,\emptyset) = 1$,
where $E_1$ and~$E_2$ are complex conjugate;
\item[(2ii)]
$W_{\Sig, \phi}(\{L-E_i-E_j,L-E_i-E_k\},0,0,\theta_1,\emptyset,\emptyset)
= 1$, where
$\{i,j,k\}\subset\{3,4,5,6\}$, the divisor $E_i$ is
real, and the divisors $E_j,E_k$ are complex conjugate; \item[(2iii)]
$W_{\Sig, \phi}(\{-(K_\Sig+E)-E_1,-(K_\Sig+E)-E_2\},0,0,\theta_1,\emptyset,\emptyset)
= 1$,
where $E_1$ and~$E_2$ are complex conjugate;
\item[(2iv)]
$W_{\Sig, \phi}(\{L-E_i-E_j,L-E_k-E_l\},0,0,\theta_1,\emptyset,\emptyset)
= -1$, where
$\{i,j,k,l\}=\{3,4,5,6\}$, $E_i, E_k$ and
$E_j, E_l$ are two pairs of complex conjugate divisors.
\item[(2v)]
$W_{\Sig,\phi}(\{-(K_\Sig+E),-(K_\Sig+E)\},0,0,\theta_2,
\emptyset,\emptyset)=1$,
if
$Q',Q''$ are complex conjugate.
\end{enumerate}
\end{enumerate}
\end{proposition}

\begin{proposition}\label{ini2}
Let $\Sig=\CS$. {\rm (}Recall that $F$ is the non-orientable
component of $\R\Sig$, and $E=L_1$.{\rm )} If $\phi$ is either $0$,
or $\phi_F$ {\rm (}{\it cf.} section \ref{W-invariants}\/{\rm )},
then among the Welschinger numbers $W_{\Sig,
\phi}(\D,\alp,\bet^{\re},\bet^{\ima},\bp^\flat,\emptyset)$, where
$(\D,\alp,\bet^\re,\bet^\ima,\bp^\flat)$ is an admissible tuple such
that
$$
\alpha,\beta^{\re}\in\Z_+^{\infty,\odd}, \;\;\;
[\D]E > 0, \,\text{and}\;\;\; R_\Sig(\D, \beta^\re + 2\beta^\ima) = 0,
$$
the only non-zero numbers are as follows:
\begin{enumerate}
\item[(1)]
in the case of a divisor class $\D = D$,
\begin{enumerate}
\item[(1i)]
$W_{\Sig, \phi}(L_i,0,\theta_1,0,\emptyset,\emptyset) = 1$, where $i=2,3$;
\item[(1ii)]
$W_{\Sig, \phi}(-(K_\Sig+E),\theta_1,\theta_1,0,\bp^\flat,\emptyset) = 1$;
\end{enumerate}
\item[(2)]
in the case of a pair $\D = (D_1, D_2)$ of divisor classes,
$$W_{\Sig, \phi}(\{L_{2i}, L_{2i+1}\},0,0,\theta_1,\emptyset,\emptyset) =
\begin{cases}
\hskip9pt 1,  & \text{if} \; \phi = \phi_F \; \text{and}
\; L_{2i} \cap L_{2i + 1} \cap F = \varnothing,\\
-1, & \text{otherwise},
\end{cases}
$$
where $i=2,3,4,5$.
\end{enumerate}
\end{proposition}

{\bf Proof of Propositions~\ref{ini1} and~\ref{ini2}}.
Both
propositions can be easily derived from Lemmas~\ref{l3}
and~\ref{l4}. We only make a comment concerning the statement of
Proposition \ref{ini2}(2). By Lemma \ref{ln1}(2), there are four
pairs $(L_{2i},L_{2i+1})$, $i=2,3,4,5$, of complex conjugate lines
crossing $E$, and the lines of each pair intersect at a real point,
which is a solitary node of the corresponding curve
$L_{2i}L_{2i+1}$. If $\phi=0$, then the contribution of such a curve
to the Welschinger number is $-1$, whereas if $\phi=\phi_F$, then by
formula (\ref{e1}) the contribution is $-1$ when the solitary node
occurs on the component $F$ and it is $1$ otherwise.  \proofend

The numbers
$W_{\Sig,\phi}(D,\alp,\bet^{\re},\bet^{\ima},\bp^\flat,\emptyset)$
in Propositions
\ref{ini1} and \ref{ini2} do not depend on the choice of $\bp^\flat$.
We skip $\bp^\flat$ and $\emptyset$ in the notation of these
numbers and simply write
$W_{\Sig,\phi}(D,\alp,\bet^{\re},\bet^{\ima})$.

\subsection{Recursive formula}\label{secRF}

\begin{theorem}\label{t1}
Let $\Sig=\PP^2_{a,b}$, $a+2b=6$, or $\CS$, let $E$ be a real
smooth $(-1)$-curve with $\R E\subset F$, and let $\phi\in
H_2(\Sigma\setminus F; \Z/2\Z)$ be a class invariant under the
action of complex conjugation, $\conj_*\phi=\phi$.

Let $D\in\Pic^\R(\Sig,E)$ be a divisor class,
and let
$\alp,\bet^{\re}\in\Z_+^{\infty,\odd}$, $\bet^{\ima}\in\Z_+^\infty$
satisfy the following conditions:
$$I(\alp+\bet^{\re}+2\bet^{\ima})=DE,\quad
R_\Sig(D,\bet^{\re}+2\bet^{\ima})>0\ .$$
Then,
$$W_{\Sig,\phi}(D,\alp,\bet^{\re},\bet^{\ima})=\sum_{k\ge 1,\
\bet^{\re}_k>0}W_{\Sig,\phi}(D,\alp+\theta_k,\bet^{\re}-
\theta_k,\bet^{\ima})$$
$$+\sum\frac{2^{\|\bet^{(0)}\|}}{\bet^{(0)}!}(l+1)
\left(\begin{matrix}\alp\\
\alp^{(0)}\alp^{(1)}...\alp^{(m)}\end{matrix}\right)
\frac{(n-1)!}{n_1!...n_m!}$$
\begin{equation}\times\prod_{i=1}^m\left(\left(
\begin{matrix}(\bet^{\re})^{(i)}\\
\gam^{(i)}\end{matrix}\right)W_{\Sig,\phi}(\D^{(i)},\alp^{(i)}
,(\bet^{\re})^{(i)},(\bet^{\ima})^{(i)})\right) \
,\label{e39}\end{equation} where $$n = R_\Sig(D,
\bet^{\re}+2\bet^{\ima}),\quad n_i=R_\Sig(\D^{(i)},
(\bet^{\re})^{(i)}+2(\bet^{\ima})^{(i)}),\ i=1,...,m\ ,
$$
and the second sum in~(\ref{e39}) is taken \begin{itemize}\item over
all integers $l\ge0$ and vectors $\alp^{(0)}\le\alp$,
$\bet^{(0)}\le\bet^{\re}$; \item over all sequences
\begin{equation}(\D^{(i)},\alp^{(i)},(\bet^{\re})^{(i)},(\bet^{\ima})^{(i)}),\
1\le i\le m \ ,\label{e40}
\end{equation}
such that, for all $i=1,...,m$,
\begin{enumerate}
\item[(1a)] $\D^{(i)}\in\Pic^\R(\Sig,E)$, and $\D^{(i)}$ is neither the
divisor class $-(K_\Sig+E)$, nor the pair
$\{-(K_\Sig + E), -(K_\Sig + E)\},$
\item[(1b)] $I(\alp^{(i)}+(\bet^{\re})^{(i)}+2(\bet^{\ima})^{(i)})=[\D^{(i)}]E$, $(\bet^{\re})^{(i)}+(\bet^{(\ima})^{(i)}\ne0$, and
$R_\Sig(\D^{(i)},(\bet^{\re})^{(i)}+2(\bet^{\ima})^{(i)})\ge0$,\end{enumerate}
and \begin{enumerate}
\item[(1c)] $D-E=\sum_{i=1}^m[\D^{(i)}]
-(2l+I\alp^{(0)}+I\bet^{(0)})(K_\Sig+E)$,
\item[(1d)] $\sum_{i=1}^m\alp^{(i)}\le\alp-\alp^{(0)}$,
\item[(1e)] $\sum_{i=1}^m(\bet^{\re})^{(i)}\ge\bet^{\re}$ and
$\sum_{i \in S}(\bet^{\ima})^{(i)}=\bet^{\ima}$,
where
$$
S=\{i\ :\ 1\le i\le m,\
\D^{(i)}\ \text{is a divisor class}\} \ ,
$$
\item[(1f)] each
tuple $(\D^{(i)},0,(\bet^{\re})^{(i)},(\bet^{\ima})^{(i)})$ with
$n_i=0$ appears in (\ref{e40}) at most once,
\end{enumerate}
\item over all sequences \begin{equation}\gam^{(i)}\in\Z_+^{\infty,
\odd},\quad \|\gam^{(i)}\|=\begin{cases}1,\ &i\in S,
\\ 0,\ &i\not\in S,\end{cases}\qquad i=1,...,m\ ,\label{e41}
\end{equation} satisfying
\begin{enumerate}\item[(2a)]
$(\bet^{\re})^{(i)}\ge\gam^{(i)}$, $i = 1$, $\ldots$, $m$, and
$\sum_{i=1}^m\left((\bet^{\re})^{(i)}-
\gam^{(i)}\right)=\bet^{\re}-\bet^{(0)}$,\end{enumerate}
\end{itemize} and the second sum in (\ref{e39}) is factorized by
simultaneous permutations in the sequences (\ref{e40}) and
(\ref{e41}).\end{theorem}

The proof is presented in section \ref{secn6}.

\begin{corollary}\label{c1} Let $\Sig=\PP^2_{a,b}$, $a+2b=6$, or $\CS$, let $E$ be a real
smooth $(-1)$-curve with $\R E\subset F$, and let $\phi\in
H_2(\Sigma\setminus F; \Z/2\Z)$ be a class invariant under the
action of complex conjugation, $\conj_*\phi=\phi$.

(1) For any divisor class $D\in\Pic^\R(\Sig,E)$ and vectors
$\alp,\bet^{\re}\in\Z_+^{\infty,\odd}$, $\bet^{\ima}\in\Z_+^\infty$
such that $I(\alp+\bet^{\re}+2\bet^{\ima})=DE$,
$R_\Sig(D,\bet^{\re}+2\bet^{\ima})\ge0$, and $\bet^{\ima}\ne 0$, one
has
\begin{equation}W_{\Sig,\phi}(D,\alp,\bet^{\re},\bet^{\ima})=0\
.\label{e43}\end{equation}

(2) For any divisor class $D\in\Pic^\R(\Sig,E)$ and vectors
$\alp,\bet\in\Z_+^{\infty,\odd}$ such that \mbox{$I(\alp+\bet)=DE$}
and $R_\Sig(D,\bet)>0$, one has
$$W_{\Sig,\phi}(D,\alp,\bet,0)=\sum_{k\ge 1,\
\bet_k>0}W_{\Sig,\phi}(D,\alp+\theta_k,\bet-\theta_k,0)$$
$$+\sum\frac{2^{\|\bet^{(0)}\|}}{\bet^{(0)}!}(l+1)\left(\begin{matrix}\alp\\
\alp^{(0)}\alp^{(1)}...\alp^{(m)}\end{matrix}\right)
\frac{(n-1)!}{n_1!...n_m!}$$
\begin{equation}\times\prod_{i=1}^m\left(\left(
\begin{matrix}(\bet^{\re})^{(i)}\\
\gam^{(i)}\end{matrix}\right)W_{\Sig,\phi}(\D^{(i)},\alp^{(i)}
,(\bet^{\re})^{(i)},(\bet^{\ima})^{(i)})\right)
,\label{e44}\end{equation} where $$n = R_\Sig(D, \bet),\quad
n_i=R_\Sig(\D^{(i)}, (\bet^{\re})^{(i)}+2(\bet^{\ima})^{(i)}),\
i=1,...,m\ ,
$$
and the second sum in~(\ref{e44}) is taken \begin{itemize}\item over
all integers $l\ge0$ and vectors $\alp^{(0)}\le\alp$,
$\bet^{(0)}\le\beta$; \item over all sequences
\begin{equation}(\D^{(i)},\alp^{(i)},(\bet^{\re})^{(i)},(\bet^{\ima})^{(i)}),\
1\le i\le m \ ,\label{e45}
\end{equation} such that, for all $i=1,...,m$,
\begin{enumerate}
\item[(1a)] $\D^{(i)}\in\Pic^\R(\Sig,E)$, and $\D^{(i)}$ is neither the
divisor class $-(K_\Sig+E)$, nor the pair
$(-(K_\Sig+E),-(K_\Sig+E))$,
\item[(1b)]
$I(\alp^{(i)}+(\bet^{\re})^{(i)}+2(\bet^{\ima})^{(i)})=[\D^{(i)}]E$,
and $R_\Sig(\D^{(i)}, (\bet^{\re})^{(i)}+2(\bet^{\ima})^{(i)})\ge0$,
\item[(1c)] $(\bet^{\ima})^{(i)}\ne0$
if and only if $\D^{(i)}$ is a pair of divisor classes, and, in such
a case, $n_i=0$, $\alp^{(i)}=(\bet^{\re})^{(i)}=0$, \end{enumerate}
and
\begin{enumerate}
\item[(1d)]
$D-E=\sum_{i=1}^m[\D^{(i)}]-
(2l+I\alp^{(0)}+I\bet^{(0)})(K_\Sig+E)$,
\item[(1e)] $\sum_{i=0}^m\alp^{(i)}\le\alp$, $\sum_{i=1}^m(\bet^{\re})^{(i)}
\ge\bet-\bet^{(0)}$,
\item[(1f)] each
tuple $(\D^{(i)},0,(\bet^{\re})^{(i)},(\bet^{\ima})^{(i)})$ with
$n_i=0$ appears in (\ref{e45}) at most once,
\end{enumerate}
\item over all sequences \begin{equation}\gam^{(i)}\in\Z_+^{\infty,
\odd},\quad
\|\gam^{(i)}\|=\begin{cases}1,\ &\D^{(i)}\ \text{is a divisor class}, \\
0,\ &\D^{(i)}\ \text{is a pair of divisor classes},\end{cases}\quad
i=1,...,m\ ,\label{e41n}
\end{equation} satisfying
\begin{enumerate}\item[(2a)] $(\bet^{\re})^{(i)}\ge\gam^{(i)}$, $i=1,...,m$, and $\sum_{i=1}^m\left((\bet^{\re})^{(i)}-
\gam^{(i)}\right)=\bet^{\re}-\bet^{(0)}$,\end{enumerate}
\end{itemize}
and the second sum in (\ref{e44}) is factorized by simultaneous
permutations in the sequences (\ref{e45}) and (\ref{e41n}).

(3) Assume that
$\phi = 0$ or $\phi = \phi_F$.
Then, all Welschinger numbers $W_{\Sig,\phi}(D,\alp,\bet,0)$, where
$D\in\Pic^\R(\Sig,E)$ is a divisor class
and $\alp, \beta \in \Z^{\infty, \odd}$
are vectors such that $I(\alp + \beta) = DE$
and $R_\Sig(D, \bet) > 0$,
are recursively determined by the formula (\ref{e44}) and
the initial conditions in Propositions \ref{ini1}, \ref{ini2}.

(4) For any divisor class $D\in\Pic^\R(\Sig,E)$, one has
\begin{equation}\WW_\phi(\Sig,D)=W_{\Sig,\phi}(D,0,(DE)\theta_1,0)\
.\label{e47}\end{equation}
\end{corollary}

{\bf Proof}. The condition
(1e)
in Theorem \ref{t1} and the
vanishing of the Welschinger numbers
$W_{\Sig,\phi}(D,\alp,\bet^{\re},\bet^{\ima})$ such that $D$ is a
divisor class, $R_\Sig(D,\bet^{\re}+2\bet^{\ima})=0$, and
$\bet^{\ima}\ne0$ (see Propositions \ref{ini1}, \ref{ini2}) directly
imply the statements (1) and (2). The claim (3) is straightforward.
Formula (\ref{e47}) comes from (\ref{e35}) and (\ref{e43}).
\proofend

\begin{remark}\label{r2}
(1) Formula (\ref{e44}) holds for any surface $\Sig'$ obtained from
$\Sig$ by successive blowing down of real $(-1)$-curves or pairs of
disjoint imaginary conjugate $(-1)$-curves: one simply has to reduce
$\Pic^\R(\Sig,E)$ to the divisor classes/pairs of divisor classes
which do not cross the blown down curves.

(2) If one blows down $L_3$ on $\CS$, the resulting surface~$\CB$ is
a conic bundle and has two real spherical components; these
components give rise to the same collections of Welschinger
invariants.
\end{remark}

\begin{corollary}\label{new-corollary}
Let $E$ be one of the three real lines of~$\CS$.

(1) For any divisor class $D\in\Pic^\R(\CS, E)$ and vectors
$\alp,\bet\in\Z_+^{\infty,\odd}$ such that \mbox{$I(\alp+\bet)=DE$}
and $R_\CS(D,\bet)>0$, one has
$$W_{\CS,\phi_F}(D,\alp,\bet,0)=\sum_{k\ge 1,\
\bet_k>0}W_{\CS,\phi_F}(D,\alp+\theta_k,\bet-\theta_k,0)$$
$$+\sum\frac{2^{\|\bet^{(0)}\|}}{\bet^{(0)}!}\left(\begin{matrix}\alp\\
\alp^{(0)}\alp^{(1)}...\alp^{(m)}\end{matrix}\right)
\frac{(n-1)!}{n_1!...n_m!}$$
\begin{equation}\times\prod_{i=1}^m\left(\left(
\begin{matrix}(\bet^{\re})^{(i)}\\
\gam^{(i)}\end{matrix}\right)
W_{\CS,\phi_F}(D^{(i)},\alp^{(i)}, \beta^{(i)}, 0)\right),\label{ecs1}
\end{equation}
where $$n = R_\CS(D, \bet),\quad
n_i=R_\CS(D^{(i)}, \beta^{(i)}),\ i=1,...,m\ ,
$$
and the second sum in~(\ref{ecs1}) is taken
\begin{itemize}
\item over
all vectors $\alp^{(0)}\le\alp$,
$\bet^{(0)} \le \bet$;
\item over all sequences
\begin{equation}
(D^{(i)},\alp^{(i)}, \bet^{(i)}, 0),\
1\le i\le m \ ,\label{ecs2}
\end{equation}
such that, for all $i=1,...,m$,
\begin{enumerate}
\item[(i)] $D^{(i)}\in\Pic^\R(\CS, E)$ is a divisor class
different from $-(K_\CS + E)$,
\item[(ii)]
$I(\alp^{(i)}+\bet^{(i)}) = D^{(i)}E$,
and $R_\CS(D^{(i)}, \bet^{(i)})\ge0$,
\end{enumerate} and
\begin{enumerate}
\item[(iii)]
$D-E=\sum_{i=1}^mD^{(i)} - (I\alp^{(0)}+I\bet^{(0)})(K_\CS + E)$,
\item[(iv)] $\sum_{i=0}^m\alp^{(i)}\le\alp$,
$\sum_{i=1}^m\bet^{(i)}\ge\bet-\bet^{(0)}$,
\item[(v)] each tuple
$(D^{(i)}, 0, \bet^{(i)}, 0)$ with $n_i=0$ appears in (\ref{ecs2})
at most once, and coincides with $(L_j, 0, \theta_1, 0)$, where $L_j
\ne E$ is a real line on~$\CS$,
\end{enumerate}
\item over all sequences
\begin{equation}
\gam^{(i)}\in\Z_+^{\infty,
\odd},\quad
\|\gam^{(i)}\| = 1, \quad
i=1,...,m\ ,\label{ecs3}
\end{equation}
satisfying
$$\bet^{(i)}\ge\gam^{(i)}, \; i = 1, \ldots, m, \;\;\;
\text{and} \;\;\; \sum_{i=1}^m\left(\bet^{(i)}-
\gam^{(i)}\right) = \bet - \bet^{(0)},$$
\end{itemize}
and the second sum in (\ref{ecs2}) is factorized by simultaneous
permutations in the sequences (\ref{ecs2}) and (\ref{ecs3}).

(2) All Welschinger numbers $W_{\CS, \phi_F}(D,\alp,\bet,0)$, where
$D\in\Pic^\R(\CS, E)$ is a divisor class
and $\alp, \beta \in \Z^{\infty, \odd}$ are vectors such that
$I(\alp + \beta) = DE$
and
$R_\CS(D, \bet) > 0$,
are recursively determined by the formula (\ref{ecs1}) and
the initial conditions in Proposition~\ref{ini2}(1).
\end{corollary}

{\bf Proof}.
Let $L_1$, $L_2$, and $L_3$ be the three real lines of $\CS$,
and let $L_j$, $j = 4$, $\ldots$,  $11$, be
the eight non-real lines of $\CS$ which intersect $L_1$
and are numbered in such a way that
$L_{2i}$ and $L_{2i+1}$ are complex conjugate
for any $i = 2, 3, 4, 5$
({\it cf}. Lemma~\ref{ln1}). Assume
that $E = L_1$.

For precisely two pairs $(L_{2i},L_{2i+1})$, $2\le i\le 5$, say, for
$i=2,3$, the intersection point belongs to $F$, and for
the other two pairs $(L_{2i},L_{2i+1})$, $i=4,5$, the intersection point
belongs to $\R\CS\setminus F$,
and hence ({\it cf}. Proposition~\ref{ini2}(2))
\begin{equation}
W_{\CS,\phi_F}(\{L_{2i},L_{2i+1}\},0,0,\theta_1)=\begin{cases}-1,\quad
& i=2,3,\\ 1,\quad & i=4,5.\end{cases}\label{eND}
\end{equation}

Combining the terms of the second sum in the righthand side of the
formula~(\ref{e44}) applied to $W_{\CS, \phi_F}(D, \alp, \beta, 0)$,
we obtain the expression
$$
\sum_{k = 0}^4\sum_{l \geq 0} x_k(l + 1)\lam_{k,l},
$$
where $k$ is the number of non-divisorial factors in a summand,
and $x_k$ is the sum of products of $k$ distinct non-divisorial terms.
Any two coefficients $\lam_{k, l}$ and $\lam_{k', l'}$ such that
$k + 2l = k' + 2l'$ coincide, so we put $\lam_{k + 2l} = \lam_{k, l}$.
One has
$$
x_0 = x_4 = 1, \;\;\; x_1 = x_3 = 0, \;\;\; x_2 = -2.
$$
Thus, the second sum in the righthand side of the formula
is equal to
$$
\sum_{l \geq 0}(l + 1)(\lam_{2l} - 2\lam_{2 + 2l} + \lam_{4 + 2l})
= \lam_0.
$$
\proofend

\begin{example}\label{ex} We present here some values of Welschinger
invariants computed by means of formulas (\ref{e44}) and
(\ref{ecs1}). In the case of $\PP^2_{a,b}$ these are the usual
Welschinger invariants; for the conic bundle $\CB$ {\rm (}see
Remark~\ref{r2}\/{\rm )} we take for~$F$ one of the components of
$\R \CB$; for $\CS$, as always, $F$ is the non-orientable component
of $\R\CS$.

\begin{center}
\begin{tabular}{|l||c|c|c|c|c|c|c|c|}
\hline $D\;\diagdown\; \Sig,\phi$ & \mbox{$\PP^2_{6,0}$} &
\mbox{$\PP^2_{4,1}$} & \mbox{$\PP^2_{2,2}$} & \mbox{$\PP^2_{0,3}$} &
\mbox{$\CB,0$} & \mbox{$\CB,\phi_F$} & \mbox{$\CS,0$}&
\mbox{$\CS,\phi_F$} \\
\hline\hline
-K & 8 & 6 & 4 & 2 & 0 & 4 & 0 & 4 \\
\hline -2K & 1000 & 522 & 236 & 78 & 0 & 512 & 0 & 160 \\
\hline
\end{tabular}
\end{center}
\end{example}

\subsection{Proof of Proposition \ref{l12} and
Theorem \ref{t1}}\label{secn6}

We simultaneously prove Proposition \ref{l12} and Theorem \ref{t1}
by induction on \mbox{$R_\Sig(D,\bet^{\re}+2\bet^{\ima})$}.

The claim of Proposition \ref{l12} for
$R_\Sig(D,\bet^{\re}+2\bet^{\ima})=0$ (the base of induction)
follows from Propositions~\ref{ini1} and~\ref{ini2}.

For induction step, we fix a tuple $(D,\alp,\bet^{\re},\bet^{\ima})$
satisfying the hypotheses of Theorem \ref{t1} and choose a divisor
class $D_0\in\Pic^\R(\Sig,E)$, $D_0\ge D$, and a sequence of arcs
$(\Lam_s)_{s=1,...,N}$ as in Lemma \ref{D0}. Pick two point
sequences $\bp^\flat=(p_s(0))_{s\in J^\flat}$ and
$\bp^\sharp=(p_s(1))_{s\in J^\sharp}$,
$J^\flat,J^\sharp\subset\{1,...,N\}$, such that
\begin{itemize} \item $s_1<s_2$ for all $s_1\in J^\flat$, $s_2\in
J^\sharp$,
\item the tuple $(D,\alp,\bet^\re,\bet^\ima,\bp^\flat)$ is admissible, \item the number of points in $\bp^\sharp$ is equal to
$R_\Sig(D,\bet^\re+2\bet^\ima)$,
\end{itemize} and prove that $W_{\Sig,\phi}(D,\alp,\bet^{\re},\bet^{\ima},\bp^\flat,\bp^\sharp)$ equals the right-hand side
of formula (\ref{e39}).

Consider the deformation diagram $\Delta$ of
$(D,\alp,\bet^\re,\bet^\ima,\bp^\flat,\widetilde\bp^\sharp,p_k(0))$
provided by~$\Lambda_k$, where $k=\min J^\sharp$ and
$\widetilde\bp^\sharp=\bp^\sharp\setminus\{p_k(1)\}$. We intend to
compute
$W_{\Sig,\phi}(D,\alp,\bet^{\re},\bet^{\ima},\bp^\flat,\bp^\sharp)$
by summing up the contributions of all connected components of
$\Delta$.

The connected components of $\Delta$ are enumerated by their roots
described in Lemma \ref{DC}.

Since $\bet^\re\in\Z_+^{\infty,\odd}$, Lemma \ref{ldediI} implies
that each connected component $\delta \subset \Delta$ with a root of
type (I) has a unique leaf, and this leaf has the same Welschinger
sign as the root. The contribution of these components of $\Delta$
gives the first summand in the right-hand side of formula
(\ref{e39}).

Let $\delta$ be a connected component of $\Delta$ with a root $C$ of
type (R). Using the description of leaves of $\delta$ given in Lemma
\ref{lDP}, we immediately conclude that the contribution of the
leaves of $\delta$ into
$W_{\Sig,\phi}(D,\alp,\bet^{\re},\bet^{\ima},\bp^\flat,\bp^\sharp)$
is as follows: \begin{itemize}\item $0$, if either there is
$(\bet^{\re})^{(i)}\not\in\Z_+^{\infty,\odd}$, or $Q',Q''$ are real
and at least one of $l',l''$ is odd; \item
$(l+1)\prod_{i=1}^mw_\phi(C^{(i)})$, if
$(\bet^{\re})^{(i)}\in\Z_+^{\infty,\odd}$ for all $i=1$, $\ldots$,
$m$, and $Q',Q''$ are complex conjugate, $l'=l''=l$;
\item
$\prod_{i=1}^mw_\phi(C^{(i)})$, if
$(\bet^{\re})^{(i)}\in\Z_+^{\infty,\odd}$ for all
$i=1$, $\ldots$, $m$, and
$Q',Q''$ are real, $l',l''$ are even.
\end{itemize}

Finally, taking into account the induction assumption of the
independence of the Welschinger numbers
$$W_{\Sig,\phi}(\D^{(i)},\alp^{(i)},(\bet^{\re})^{(i)},
(\bet^{\ima})^{(i)},(\bp^\flat)^{(i)},(\bp^\sharp)^{(i)}),
\quad i=1, \ldots, m\ ,$$
on the choice of pairs of point sequences
$(\bp^\flat)^{(i)},(\bp^\sharp)^{(i)}$ in $D_0$-CH position,
$i=1$, $\ldots$, $m$,
and summing up over connected components of $\Delta$,
we immediately obtain formula (\ref{e39}) and, in virtue of this
formula, also the independence of
$W_{\Sig,\phi}(D,\alp,\bet^{\re},\bet^{\ima},\bp^\flat,\bp^\sharp)$
on the choice of $\bp^\flat,\bp^\sharp$.\proofend

\section{Applications}\label{applications}

\subsection{Positivity and asymptotics}\label{sec-pa}
A divisor class $D$ on a surface $\Sig$ is called {\it nef}\; if $D$
non-negatively intersects any algebraic curve on $\Sig$. A nef
divisor class $D$ is big if $D^2>0$.

\subsubsection{The case of $\Sig=\PP^2_{a,b}$ with $a+2b\le6$, $b\le2$} 
In this case, the real part $\R\Sig$ is nonempty and connected, and
hence we can speak only of the usual Welschinger invariants, which
we simply denote by $\WW(\Sig,D)$ omitting $\phi$ in the notation.

\begin{theorem}\label{tn5}
Let $\Sig=\PP^2_{a,b}$, where $a+2b\le 6$, $b\le2$.
Then, for any
nef and big real divisor class $D$ on $\Sig$,
\begin{itemize}\item
the invariant $\WW(\Sig,D)$ is positive; in particular, through any
generic collection of $c_1(\Sig)D-1$ real points in $\Sig$ one can
trace a real rational curve $C\in|D|$, \item the following
asymptotic relation holds:
\begin{equation}\
\log\WW(\Sig,nD)=c_1(\Sig)D\cdot n\log n+O(n),\quad n\to +\infty\
,\label{en12}
\end{equation}
and in particular,
\begin{equation}
\lim_{n\to+\infty}\frac{\log W(\Sig,nD)}
{\log GW(\Sig,nD)}=1\ .\label{en12a}
\end{equation}
\end{itemize}
\end{theorem}

\begin{remark}\label{rn2}

(1) Theorem \ref{tn5} covers all the cases studied
in \cite{IKS2,IKS4,IKS6,Sh2}.
The statement holds true for the surface
$\PP^2_{0,3}$ as well, but the proof requires another
approach and will be presented in a forthcoming paper.

(2) If $D$ is not nef or not big, then $\WW(\Sig,D)=1$ or $0$
depending on whether or not the linear system $|D|$ contains an irreducible
curve (for the existence of rational irreducible
representatives see, for instance, \cite{GLS}).

(3) The positivity statement and the existence of real rational
curves do not extend in the same form to all real unnodal
Del Pezzo
surfaces, for example, $\WW(\PP^2_{0,4},-K_{\PP^2_{0,4}})=0$, and
there are generic
configurations of four pairs of imaginary conjugate blown
up points such that the linear system $|-K_{\PP^2_{0,4}}|$ does not
contain any real rational curve (cf. \cite[Section 3.1(1)]{IKS2}).
\end{remark}

\begin{lemma}\label{new-lemma}
Let $\Sig$ be an unnodal Del Pezzo surface,
and $D \in \Pic(\Sig)$.
\begin{enumerate}
\item[(i)]
The divisor class $D$ is nef and nonzero if and
only if its intersection with any $(-1)$-curve
on~$\Sig$ is non-negative and
$K_\Sig D<0$.
\item[(ii)]
Assume that $\Sig = \PP^2_{a, b}$, where $a + 2b \leq 6$.
If $D$ is nef and nonzero, then $D^2\ge 0$, and $D$ can be
represented by a union of rational curves, which are real if $D$ is
real.
More precisely, if $D^2>0$, then~$D$ can be represented by an
irreducible rational curve; if $D^2=0$, then $D=kD'$, where $D'$ is
primitive {\rm (}not divisible by a natural number $>1$\/{\rm )},
$D'$ can be
represented by an irreducible nonsingular rational curve, and $D$ can be
represented by a union of $k$ disjoint irreducible nonsingular rational
curves.
\end{enumerate}
\end{lemma}

{\bf Proof}.
The statements follow, for instance, from
\cite[Theorems 5.1, 5.2, and Remark 5.3]{GLS}.
In particular, if~$D$ is real,
the construction in
\cite[Section 5.2]{GLS}
gives real representatives.
\proofend

{\bf Proof of Theorem~\ref{tn5}}.
The case of  $\Sigma=\PP^2$
is settled in
\cite{IKS,IKS2}, so we assume that $a+b>0$.

Applying 
induction on $-K_\Sig D$, we prove that
\begin{equation}
\WW(\Sig,D) 
>0\label{eP}\end{equation}
for all nef and big real divisor classes $D$ on $\Sig$. The base of induction, that is
the case
$-K_\Sig D=1$, is provided by
Proposition \ref{ini1}. Suppose that $-K_\Sig D>1$ and perform the induction step.

We claim that either there exists a pair of disjoint complex conjugate
$(-1)$-curves, which do not intersect $D$, or $D$ intersects with all imaginary $(-1)$-curves
and there exists a real $(-1)$-curve $E$ such that $DE'>0$ for each real $(-1)$-curve
$E'$
intersecting $E$.
Indeed, if $DE'=0$ for a non-real $(-1)$-curve
$E'$, then $D\overline{E'}=0$,
and we
have $E'\overline{E'}=0$
(otherwise, one would have 
$(E'+\overline{E'})^2\ge 0$ yielding $D^2\le0$). Since $a+b>0$, there are
real $(-1)$-curves, and
the same argument shows that if $DE=0$ for some real $(-1)$-curve $E$, then
$DE'>0$ for each real $(-1)$-curve
$E'$
that intersects $E$.

Now, in the presence of a complex conjugate pair $E',\overline{E'}$ of disjoint
$(-1)$-curves such that $DE'=D\overline{E'}=0$, we blow down both $E'$ and
$\overline{E'}$ and obtain
a surface of one of the types $\PP^2_{a',b'}$, $a'+2b'\le4$, $b'\le1$,
for which the relations (\ref{eP}) and (\ref{en12}) follow from
\cite{IKS,IKS2,IKS4,IKS3,IKS6,Sh2}.

Thus, suppose that $D$ intersects with each imaginary $(-1)$-curve, and
that there exists a real $(-1)$-curve $E$ such that $DE'>0$ for all real
$(-1)$-curves $E'$ intersecting $E$.
Notice that
$$\WW(\Sig,D) = W_\Sig(D,0,(DE)\theta_1,0)$$ and
apply formula (\ref{e44}).
By Proposition \ref{ini1},
all the initial Welschinger numbers are non-negative, and hence, due
to the positivity of the coefficients in formula (\ref{e44}), we get
that $$W_\Sig(D',\alp,\bet,0)\ge 0\quad\text{for any}\quad
D'\in\Pic^\R(\Sig,E),\ \alp,\bet\in\Z^{\infty,\odd}_+,\
I(\alp+\bet)=D'E\ .$$ Putting
$W_\Sig(D,(DE+1)\theta_1,-\theta_1,0)=0$ and applying $(DE + 1)$
times the formula (\ref{e44}), we get
$$W_\Sig(D,k\theta_1,(DE-k)\theta_1,0)=W_\Sig(D,(k+1)\theta_1,
(DE-k-1)\theta_1,0) + \tau_k,\quad k=0, \ldots, DE\ ,$$ where
$\tau_k$ stands for the second sum in the right-hand side
of~(\ref{e44}). One has
\begin{equation}
W_\Sig(D,0,(DE)\theta_1,0) = \tau_0 + \ldots + \tau_{DE}\ ,\label{eS}
\end{equation}
and to prove the positivity of $W_\Sig(D, 0, (DE)\theta_1, 0)$, it
is sufficient to find at least one positive term in $\tau_0 + \ldots
+ \tau_{DE}$.

Since $E$
(as any $(-1)$-curve on a del Pezzo surface of degree $\ge 3$)
crosses any other $(-1)$-curve in at most one point,
the divisor class $D-E$ non-negatively intersects
each $(-1)$-curve,
and, thus, $D-E$
is nef
and nonzero (see Lemma~\ref{new-lemma}(i)).

In addition, $-K_\Sig(D-E)<-K_\Sig D$
and
$$\tau_0 \ge (DE+1)W_\Sig(D-E,0,(DE+1)\theta_1,0).$$
Hence, if $D-E$ is big, we can replace $D$ with $D - E$ in our
procedure.

Otherwise, according to Lemma~\ref{new-lemma}, one has $D-E=kD'$
with a nef primitive $D'$ such that $(D')^2=0$ and $D'E>0$. Since
$R_\Sig(D',(D'E)\theta_1)=\dim|D'|=1$, we have
$R_\Sig(D',(D'E-1)\theta_1)=0$, and then from Proposition
\ref{ini1}, we can derive that
$W_\Sig(D',\theta_1,(D'E-1)\theta_1,0)=1$. Formula (\ref{e44}) gives
then
$$W_\Sig(D',0,(D'E)\theta_1,0)\ge
W_\Sig(D',\theta_1,(D'E-1)\theta_1,0)>0\ .$$ This implies that the
term $\tau_{k - 1}$ in (\ref{eS}) is positive, since this term
contains the summand $(W_\Sig(D',0,(D'E)\theta_1,0))^k>0$; hence
(\ref{eP}).

\smallskip
Due to the upper bound $W(\Sig,D)\le
GW(\Sig,D)$ and the asymptotics \mbox{$\log
GW(\Sig,nD)=c_1(\Sig)D\cdot n\log n+O(n)$} (see~\cite{IKS5}), to
prove the asymptotic relations~(\ref{en12})
and~(\ref{en12a}), it is
sufficient
to establish for any nef and big real divisor class $D$
the inequality
\begin{equation}\log W_\Sig(nD,0,n(DE)\theta_1,0)\ge
c_1(\Sig)D\cdot n\log n+O(n),\quad n\to +\infty\ .\label{e48}
\end{equation}
As we noticed above, we can suppose that $D$ intersects with each imaginary $(-1)$-curve, and
that there exists a real $(-1)$-curve $E$ such that $DE'>0$ for all real
$(-1)$-curves $E'$ intersecting $E$. We claim that
there exist positive numbers $\xi,\eta,\zeta$
(depending on $D$) such that
\begin{equation}W(\Sig,nD)=W_\Sig(nD,0,n(DE)\theta_1,0)\ge
\xi\eta^{n\zeta}(-nK_\Sig D-1)! \ ,\label{eA}\end{equation} for any
$n \ge 1$, which clearly implies (\ref{e48}). Indeed, $D-E$ is nef.
If $D-E$ is
big, then, applying formula (\ref{eS}) to
$W_\Sig((n+1)D,0,(n+1)(DE)\theta_1,0)$, we obtain that the term
$\tau_2$ contains the sum
$$\frac{(-(n+1)K_\Sig
D-2)!}{(-K_\Sig(D-E)-1)!}(DE+1)\cdot
W_\Sig(D-E,0,(DE+1)\theta_1,0)$$
$$\times\sum_{i=1}^{n-1}\frac{i(n-i)(DE)^2\cdot W(iD,0,i(DE)\theta_1,0)
\cdot W_\Sig((n-i)D,0,(n-i)(DE)\theta_1,0)}{(-iK_\Sig
D-1)!(-(n-i)K_\Sig D-1)!}\ ,$$
which means that
the sequence
$$u_n = \frac{nW_\Sig(nD, 0, n(DE)\theta_1, 0)}{(-nK_\Sig D - 1)!}$$
satisfies the condition
$$
u_n \geq \sum_{i = 1}^{n - 1} c u_i u_{n - i}
$$
for a certain positive constant~$c$. As is well known ({\it
cf}.~\cite{Itz}), this implies the inequality~(\ref{eA}).

If $D - E$ is not big, then, according to Lemma~\ref{new-lemma},
one has
$D-E=kD'$, where $D'$ is a nef divisor class
with $D'E>0$, $(D')^2=0$, and
$W_\Sig(D',0,(D'E)\theta_1,0)=1$. Again applying formula (\ref{e44})
to $W_\Sig((n+1)D,0,(n+1)(DE)\theta_1,0)$, we obtain that the term
$\tau_{k + 1}$ contains the sum $$(-(n+1)K_\Sig D-2)!(D'E)^k$$
$$\times\sum_{i=1}^{n-1}\frac{i(n-i)(DE)^2\cdot W(iD,0,i(DE)\theta_1,0)
\cdot W_\Sig((n-i)D,0,(n-i)(DE)\theta_1,0)}{(-iK_\Sig
D-1)!(-(n-i)K_\Sig D-1)!}\ ,$$ which, as above, implies the
inequality~(\ref{eA}). \proofend

\subsubsection{The case of $\Sig=\CS$}\label{sec-sss3}
Recall that $F$ denotes the non-orientable component of $\R\CS$.

\begin{theorem}\label{t5} For any nef and big real divisor class $D$
on $\CS$, the Welschinger invariant $\WW_{F}(\CS,D)$ is
positive. In particular, through any collection of $c_1(\CS)D-1$
generic points of $F$ one can trace a real rational curve $C\in|D|$
passing through the given points. Furthermore, one has $$\log
W_{F}(\CS,nD)=c_1(\CS)D\cdot n\log n+O(n),\quad n\to \infty\
.$$ In particular,
$$\lim_{n\to\infty}\frac{\log\WW_{F}(\CS,nD)}{\log GW(\PP^2_6,nD)}=1\
.$$
\end{theorem}

\begin{remark}\label{rsss3}
(1) The positivity of the usual Welschinger invariants $\WW(\CS, D)$
does not hold: $\WW(\CS,-K_{\CS})=0$, whereas
$\WW_F(\CS,-K_{\CS})=4$ {\rm (}{\it cf.} Example \ref{ex}\/{\rm )}.
Indeed, for any $2 = K^2_\CS - 1$ generic points of $F$, there are
exactly two planes passing through these points and tangent to the
spherical component of $\R\CS$; each of these two planes intersects
$\CS$ along a cubic with a solitary node belonging to the spherical
component. The planes passing through the chosen two points and
tangent to~$F$ give rise to rational cubics whose total contribution
to $\WW(\CS, -K_{\CS})$ {\rm (}as well as to $\WW_F(\CS,
-K_{\CS})${\rm )} is equal to the Euler characteristics of~$F$ blown
up at $3$ points, wherefrom the above total values follow.

(2) For any real Del Pezzo surface $\Sig$ with disconnected real
part and for any connected component~$F$ of $\R\Sig$, one has
$W(\Sig, -K_\Sig, F) \leq 0$. Indeed, $W(\Sig, -K_\Sig, F) =
-\chi(\R\Sig) + K^2_\Sig$, where $\chi(\R\Sig)$ is the Euler
characteristics of $\R\Sig$, and the inequality $-\chi(\R\Sig) +
K^2_\Sig\leq 0$ follows, for example, from Comessatti's
classification of real rational surfaces~\cite{Com} (see
also~\cite{DK}).

(3) Theorem \ref{t5} implies a similar statement for the invariants
$W_F(\CB,D)$, where $F$ is any of the two connected components of
$\R \CB$ (see Remark \ref{r2}): the nef and big real divisor classes
on $\CB$ can be viewed as the divisor classes
$D=d_1L_1+d_2L_2+d_3L_3$ on $\CS$, where $L_1$, $L_2$, and $L_3$ are
the three real lines on $\R\CS$ and $d_1=d_2+d_3$, $d_2,d_3>0$.
\end{remark}

{\bf Proof of Theorem~\ref{t5}}.
Recall that a real
divisor class $D$ on $\CS$ is nef and big if and only if
$D=d_1L_1+d_2L_2+d_3L_3$, where
$$d_1,d_2,d_3>0\quad\text{and}\quad d_i+d_j\ge d_k,\
\{i,j,k\}=\{1,2,3\}\ .$$

Let $D = d_1L_1 + d_2L_2 + d_3L_3$, and
$d_i=\max\{d_1,d_2,d_3\}>1$. Applying formula (\ref{ecs1}) with
$E=L_i$, we
obtain
$$\displaylines{
W_F(\CS,D)=W_{\CS,\phi_F}(D,0,(d_1+d_2+d_3-2d_i)\theta_1,0) \cr \ge(d_1+d_2+d_3-2d_i+1)W_{\CS,\phi_F}
(D-E,0,(d_1+d_2+d_3-2d_i+1)\theta_1,0) \cr
\geq
W_F(\CS,D-E)
}
$$
with a nef and big real divisor class $D-E$. After
$d_1+d_2+d_3-3$ steps we get
$$W_F(\CS,D)\ge W_F(\CS,L_1+L_2+L_3)\ .$$
Since $W_F(\CS, L_1 + L_2 + L_3) = 4$ (see Remark~\ref{rsss3}(1)),
one has $W_F(\CS, D) > 0$ for any nef and big real divisor~$D$
on $\CS$.

Assume now that $d_1\ge d_2\ge d_3>0$, and pick a number $n\in\N$.
Putting $E = L_1$ and applying $n(d_1-d_2)$ times the
formula~(\ref{ecs1}), we get
$$W_F(\CS,nD)\ge
\frac{(nd_3)!}{(n(d_2+d_3-d_1))!}W_F(\CS,nD'),
$$
where $D' = d_2(L_1 + L_2) + d_3L_3$.
Then, putting alternatively $E = L_1$
and $E = L_2$, we apply
$2n(d_2-d_3)$ times the formula~(\ref{ecs1}) and get
$$W_F(\CS,nD')\ge
\left(nd_3(nd_3+1)\right)^{n(d_2-d_3)}W_F(\CS, nD''),$$
where $D'' = d_3(L_1 + L_2 + L_3)$.
Finally, putting cyclically $E = L_1$, $E = L_2$, and $E = L_3$,
we apply $3nd_3-3$ times
the formula~(\ref{ecs1}) and get
$$W_F(\CS,nD'')\ge
\frac{(nd_3+1)!(nd_3)!(nd_3-1)!}{2}W_F(\CS, L_1 + L_2 + L_3).$$
The above inequalities give
$$\log
W_F(\CS,nD)\ge(d_1+d_2+d_3)n\log n+O(n)=c_1(\CS)D\cdot n\log
n+O(n)$$ which implies the desired asymptotics. \proofend

\subsection{Monotonicity}\label{sec-mon}

\begin{lemma}\label{ln2} (1)
Let $\Sig$ be an unnodal Del Pezzo
surface, and $D$, $D'$ be nef and big divisor classes on $\Sig$.
If $D - D'$ is effective, then $D - D'$ can be decomposed into a
sum $E^{(1)}+...+E^{(k)}$ of smooth rational $(-1)$-curves such that
each of $D^{(i)}=D' +\sum_{j\le i}E^{(j)}$ is nef and big, and
satisfies $D^{(i)}E^{(i+1)}>0$, $i=0, \ldots,k-1$.

(2) Let $D$ and $D'$ be nef and big real divisor classes on $\CS$.
If $D - D'$ is effective, then $D - D' = E^{(1)} + \ldots +
E^{(k)}$, where $E^{(1)}, \ldots, E^{(k)}\in\{L_1,L_2,L_3\}$, and
the following properties hold for any $i = 0$, $\ldots$, $k - 1$:
the divisor class $D^{(i)}=D' + \sum_{j\le i}E^{(j)}$ is nef and
big, and $D^{(i)}E^{(i+1)}>0$.
\end{lemma}

{\bf Proof}. The proof of the first claim literally coincides with
the proof of \cite[Lemma 30]{IKS6}. To prove the second statement,
write $D = d_1L_1 + d_2L_2 + d_3L_3$ and $D' = d'_1L_1 + d'_2L_2 +
d'_3L_3$, where $d_j$, $j = 1, 2, 3$, and $d'_j$, $j = 1, 2, 3$, are
positive integers such that
$$\displaylines{
d_1 + d_2 \geq d_3, \;\;\; d_1 + d_3 \geq d_2,
\;\;\; d_2 + d_3 \geq d_1, \cr
d'_1 + d'_2 \geq d'_3, \;\;\; d'_1 + d'_3 \geq d'_2,
\;\;\; d'_2 + d'_3 \geq d'_1.
}
$$
Put $k = (d_1 + d_2 + d_3) - (d'_1 + d'_2 + d'_3)$ and $D^{(k)} =
D$. Define inductively $E^{(i + 1)}$ and $D^{(i)}$, $i = k - 1$,
$\ldots$, $0$, in such a way that each $E^{(i + 1)}$ is a real line
of~$\CS$, each divisor class $D^{(i)}$ is nef and big, $D^{(0)} =
D'$, and $D^{(i)}E^{(i + 1)} > 0$, for any $i = k - 1$, $\ldots$,
$0$. This can be done as follows. Write $D^{(i + 1)}$ in the form
$d^{(i+1)}_1L_1 + d^{(i+1)}_2L_2 + d^{(i+1)}_3L_3$, where
$d^{(i+1)}_1$, $d^{(i+1)}_2$, and $d^{(i+1)}_3$ are positive
integers, and choose among the coefficients $d^{(i+1)}_j$ such that
$d^{(i+1)}_j > d'_j$ a maximal one, $d^{(i+1)}_{j^{(i+1)}}$. Put
$E^{(i+1)} = L_{j^{(i+1)}}$ and $D^{(i)} = D^{(i+1)} - E^{(i+1)}$.
\proofend

\begin{theorem}\label{tn7}
(1) Let $D$ and $D'$ be nef and big divisor classes on $\PP^2_{6,0}$
such that $D - D'$ is effective. Then $W(\PP^2_{6,0}, D)\ge
W(\PP^2_{6,0}, D')$. Moreover, in the notation of
Lemma~\ref{ln2}(1), one has
$$W(\PP^2_{6,0}, D) \ge \prod_{i=1}^k(D^{(i-1)}E^{(i)}) \cdot W(\PP^2_{6,0}, D')\ .$$

(2) Let $D$ and $D_2$ be nef and big real divisor classes on $\CS$
such that $D - D'$ is effective.
Then $W_F(\CS, D) \ge
W_F(\CS, D')$.
Moreover, in the notation of Lemma~\ref{ln2}(2),
one has
$$W_F(\CS, D) \ge \prod_{i=1}^k(D^{(i-1)}E^{(i)}) \cdot W_F(\CS, D')\ .$$
\end{theorem}

{\bf Proof}.
The statements immediately follow from
formulas (\ref{e44}) and (\ref{ecs1}). \proofend

\subsection{Mikhalkin's congruence}\label{secn22}

\begin{theorem}\label{t8}
Let $\Sig = \PP^2_{6, 0}$. Then, for any nef and big divisor class
$D$ on $\Sig$, one has
\begin{equation}
W(\Sig, D) = GW(\Sig, D) \mod 4 \ .\label{congr}
\end{equation}
\end{theorem}

{\bf Proof}. Let~$E$ be a $(-1)$-curve on~$\Sig$. For any big and
nef divisor class $D$ on $\Sig$ and any sequences $\alp, \beta \in
\Z^\infty_+$ such that $I(\alp + \beta) = DE$, consider a generic
collection $\bz^\flat=\{z_{i,j}\ :\ i\ge1,\ 1\le j\le\alp_i\}$ of
$\|\alp\|$ points on $E$, and the variety
$V^\C_\Sig(D,\alp,\bet,\bz^\flat)$ which is the union of
$R_\Sig(D,\bet)$-dimensional components of the family of complex
reduced irreducible rational curves $C\in|D|$ which intersect $E$ in
the following way:
\begin{itemize}
\item $C$ has one local branch
at each of the points of $C\cap E$,
\item $(C\cdot E)_{z_{i, j}} = i$ for all
$i \ge 1$, $1\le j \le \alp_i$,
\item for each $i \ge 1$, there are
precisely $\bet_i$ points $q\in(C\cap E)\setminus\bz^\flat$ such
that $(C\cdot E)_q = i$
\end{itemize}
({\it cf.} Section~\ref{Welsch-inv} and \cite[Definition 2.4]{MS}).
Denote by $N_\Sig(D, \alp, \beta)$ the degree of $V^\C_\Sig(D, \alp,
\beta, \bz^\flat)$. In particular, $N_\Sig(D,0,(DE)\theta_1) =
GW(\Sig, D)$ for any nef and big divisor class $D$.

We prove the following statement:
\begin{equation}
W_\Sig(D,\alp,\bet,0) = I^\bet N_\Sig(D,\alp,\bet) \mod 4\
,\label{e112}\end{equation} for any divisor class
$D\in\Pic^\R(\Sig,E)$ and any sequences
$\alp,\bet\in\Z^{\infty,\odd}_+$ such that $I(\alp+\bet)=DE$. This
statement immediately implies the statement of the theorem.

Using induction on $R_\Sig(D,\bet)$ and the recursive formula
(64) from \cite{MS}, we easily derive that the numbers
$N_\Sig(D,\alp,\bet)$ are even if
$\bet\not\in\Z^{\infty,\odd}_+$, and hence
\begin{equation}I^{\bet}\cdot N_\Sig(D,\alp,\bet)
= 0\mod4\quad\text{if}\ \bet\not\in
\Z_+^{\infty,\odd}\ .\label{eCO}
\end{equation}
Then, using Proposition \ref{ini1}, we check the congruence
(\ref{e112}) in the case $R_\Sig(D,\bet)=0$. Finally, we proceed by
induction on $R_\Sig(D,\bet)$, comparing term by term \cite[Formula
(64)]{MS} and formula (\ref{e44}) and using the following
observations:
\begin{itemize}
\item $\Pic^\R(\Sig, E)$ contains only divisor classes,
and hence the parameters $(\bet^{\ima})^{(i)}$ in (\ref{e44}) always
vanish,
\item for any integer~$j$, one has
$$j^2 = \begin{cases}0\mod4,\
& \text{if} \; j = 0\mod2,\\ 1\mod4,\
& \text{if} \; j = 1 \mod 2,\end{cases}$$
\item for any non-negative integer~$k$, one has
$$\binom{k+3}{3} = \begin{cases}0\mod4,\ &\text{if}\
k = 1\mod2,\\ l+1\mod4,\
&\text{if}\ k=2l,\end{cases}$$
\item for any sequence $\beta^{(0)} \in \Z^{\infty, \odd}_+$, one has
$$2^{\|\bet^{(0)}\|}I^{\bet^{(0)}} = 2^{\|\bet^{(0)}\|} \mod 4,
$$
\end{itemize}
\proofend

The congruence (\ref{congr}) was established by G. Mikhalkin
(\cite{Mip}, {\it cf}.~\cite{BM}) for
$\Sig = \PP^2$,
$\PP^1 \times \PP^1$, and $\PP^2_{a,0}$, $a=1,2,3$,
and then extended to the cases of $\PP^2_{4, 0}$ and $\PP^2_{5, 0}$
in~\cite{IKS6}.

{\ncsc
Universit\'e Pierre et Marie Curie and
Institut Universitaire de France\\[-21pt]

Institut de Math\'ematiques de Jussieu\\[-21pt]

4 place Jussieu,
75005 Paris, France}\\[-21pt]

{\it E-mail}:
{\ntt itenberg@math.jussieu.fr}

\vskip10pt

{\ncsc Universit\'{e} de Strasbourg and IRMA \\[-21pt]

7, rue Ren\'{e} Descartes, 67084 Strasbourg Cedex, France} \\[-21pt]

{\it E-mail}: {\ntt viatcheslav.kharlamov@math.unistra.fr}

\vskip10pt

{\ncsc School of Mathematical Sciences \\[-21pt]

Raymond and Beverly Sackler Faculty of Exact Sciences\\[-21pt]

Tel Aviv University,
Ramat Aviv, 69978 Tel Aviv, Israel} \\[-21pt]

{\it E-mail}: {\ntt shustin@post.tau.ac.il}


\begin{thebibliography}{99}

\bibitem{Bru} Brugall\'e, E.,
private communication.

\bibitem{BM} Brugall\'e, E., Mikhalkin, G.:
Enumeration of curves via floor diagrams. C. R. Math. Acad. Sci.,
Paris {\bf 345}, no. 6, (2007), 329--334.

\bibitem{CH} Caporaso, L., and Harris, J.:
Counting plane curves of any genus.
{\it Invent. Math.} {\bf 131} (1998), no. 2, 345--392.

\bibitem{Com} Comessatti, A;:
Fondamenti per la geometria sopra le superficie razionali dal punto di
vista reale.
{\it Math. Ann.} {\bf 73} (1912), 1--72.

\bibitem{DK} Degtyarev, A., Kharlamov, V.:
Around real Enriques surfaces. {\it Revista Matematica de la
Universidad Complutense de Madrid} {\bf 10} (1997), 93--109.

\bibitem{Itz}
Di Francesco, P., Itzykson, C.:
Quantum intersection rings.
{\it in ``The moduli space of curves''}
(Texel Island, 1994),  81--148, Progr. Math., 129,
Birkh\"{a}user Boston, Boston, MA, 1995.

\bibitem{GLS} Greuel, G.-M., Lossen, C., and Shustin, E.: Geometry of families
of nodal curves on the blown up projective plane. {\it Trans. Amer.
Math. Soc.} {\bf 350} (1998), no. 1, 251--274.

\bibitem{IKS} Itenberg, I., Kharlamov, V., and Shustin, E.: Welschinger
invariant and enumeration of real rational curves. {\it
International Math. Research Notices} {\bf 49} (2003), 2639--2653.

\bibitem{IKS2} Itenberg, I.,
Kharlamov, V., and Shustin, E.: Logarithmic equivalence of
Welschinger and Gromov-Witten invariants. {\it Russian Math.
Surveys} {\bf 59} (2004), no. 6, 1093--1116.

\bibitem{IKS5} Itenberg, I., Kharlamov, V., and Shustin, E.: Logarithmic
asymptotics of the genus zero Gromov-Witten invariants of the blown
up plane. {\it Geometry and Topology} {\bf 9} (2005), paper no. 14,
483--491.

\bibitem{IKS4}  Itenberg, I.,
Kharlamov, V., and Shustin, E.: New cases of logarithmic equivalence
of Welschinger and Gromov-Witten invariants. {\it Proc. Steklov
Math. Inst.} {\bf 258} (2007), 65--73.

\bibitem{IKS3} Itenberg, I.,
Kharlamov, V., and Shustin, E.: A Caporaso-Harris type formula for
Welschinger invariants of real toric Del Pezzo surfaces. {\it
Comment. Math. Helv.} {\bf 84} (2009), 87--126.

\bibitem{IKS6} Itenberg, I.,
Kharlamov, V., and Shustin, E.: Welschinger invariants of small
non-toric Del Pezzo surfaces. {\it J. of the European Math. Soc.}
{\bf 15} (2013), no. 2, 539--594.

\bibitem{Mip} Mikhalkin,~G.: Private communication.

\bibitem{Seg} Segre,~B.:
{\it The nonsingular cubic surfaces}. Clarendon Press, Oxford, 1942.

\bibitem{MS} Shoval,~M., and Shustin,~E.: On Gromov-Witten invariants of del
Pezzo surfaces. {\it Int. J. Math.} {\bf 24} (2013), no.7, 44 pp. DOI:
10.1142/S0129167X13500547. 

\bibitem{Sh2} Shustin,~E.: Welschinger invariants of toric Del Pezzo
surfaces with non-standard real structures. {\it Proc. Steklov Math.
Inst.} {\bf 258} (2007), 218--247.

\bibitem{Va}
Vakil, R.: Counting curves on rational surfaces. {\it Manuscripta
Math.} {\bf 102} (2000), 53--84.

\bibitem{W1} Welschinger,~J.-Y.:
Invariants of real symplectic 4-manifolds and lower bounds in real
enumerative geometry. {\it Invent. Math.} {\bf 162} (2005), no. 1,
195--234.

\end{thebibliography}
\end{document}